\documentclass[draft]{article}
\def\today{14.6.10} 
\usepackage{amsmath,amsfonts,amsthm,amssymb,amscd}

\renewcommand{\Re}{\mathop{\rm Re}\nolimits}
\renewcommand{\Im}{\mathop{\rm Im}\nolimits}
\def\S{\mathhexbox278}

\newcommand{\beq}{\begin{equation}}
\newcommand{\ee}{\end{equation}}

\theoremstyle{plain} \newtheorem{theorem}{Theorem}[section]
\newtheorem{lemma}[theorem]{Lemma}
\newtheorem{proposition}[theorem]{Proposition}
 \theoremstyle{definition}
\newtheorem{definition}[theorem]{Definition} \theoremstyle{remark}
\newtheorem{remark}[theorem]{Remark}

\newcommand{\R}{{\mathbb R}} \newcommand{\U}{{\mathcal U}}

\newcommand{\Z}{{\mathbb Z}}

\newcommand{\Tr}{{\mathcal T}}
\newcommand{\Tra}{{\phi}}

\newcommand{\resto}{{\mathcal R}}
\def\im{{\rm i}}

\newcommand{\C}{\mathbb{C}}

\font\strana=cmti10
\def\lie{\hbox{\strana \char'44}}
\def\uno{{\kern+.3em {\rm 1} \kern -.22em {\rm l}}}

\def\norma#1{\left\| #1\right\|}

\numberwithin{equation}{section}

\begin{document}

\title{On scattering of small energy solutions of non autonomous
hamiltonian nonlinear Schr\"odinger equations}

\author {Scipio Cuccagna}

\date{\today}
\maketitle
\begin{abstract}  We revisit a result by Cuccagna, Kirr and Pelinovsky about   the
cubic nonlinear Schr\" odinger equation (NLS) with an attractive
localized potential and a time-dependent factor in the nonlinearity.
We show
  that, under generic  hypotheses on the linearization  at 0 of
 the equation,    small energy solutions are asymptotically free. This is yet a
     new
     application of the hamiltonian structure, continuing a program
     initiated in a paper by Bambusi and Cuccagna.
\end{abstract}

\section{Introduction}
\label{section:introduction} We consider for $ (t,x)\in\mathbb{
R}\times
 \mathbb{ R}^3$ the   nonlinear
Schr\"odinger equation  \eqref{NLS1}

\begin{equation}\label{NLS1}
 \im   u _t(t,x)= \mathcal{H}u (t,x)+\gamma(t)|u (t,x)|^2u(t,x) ,
 \, u (0,x)=u_0(x).
\end{equation}
Here $\mathcal{H}:= -\Delta +V(x)+\underline{c}$ with
$\underline{c}>0$   a constant. $\gamma(t)$ is of form
\begin{equation} \label{eq:nonlinearity} \gamma(t) =
\gamma_0+\gamma_1\cos(  t), \qquad \gamma_0,\ \gamma_1\in
\mathbb{R}, \, \gamma_1\neq 0   .\end{equation} We will assume the
following hypotheses.

\begin{itemize}
\item[(H1)]  $V(x)$ is a real valued Schwartz
function.

\item[(H2)]   We assume $\mathcal{H}
\ge 0$.

\item[(H3)] The set of eigenvalues
$\sigma _d(\mathcal{H}) $ is  contained in $[0 ,\underline{c})$.
Specifically, we assume that $0\in \sigma _d(\mathcal{H}) $ and
that the sum of the multiplicities of the eigenvalues is $n+1$. We
write $ 0=\lambda _0 <\lambda _1\le  \dots \le \lambda _n $
repeating each eigenvalue a number of times equal to its
multiplicity (notice that it is well known  that
 (H1)--(H2) implies $\dim  \ker\mathcal{H} \le 1$).

\item[(H4)]  $\underline{c}$ is not an eigenvalue or a
resonance for   $\mathcal{H}$, i.e. there are no   nonzero solutions
of $\Delta u=Vu$ in $ \mathbb{R}^3$ with $|u(x)|\lesssim \langle
x\rangle ^{-1}.$

\item[(H5)] $\underline{c}\not \in  \mathbb{N}$.

\item[(H6)]   $\forall$ $j=1,..., n$   there exists    $N_j\in \mathbb{N}$
 such that $  N_j \lambda _j<\underline{c}<  (N_j+1)\lambda _j$. Notice that  $N_1=\sup _jN_j$.

\end{itemize}

Let now $[\underline{c}]\in \Z$ be the integral part of $\underline{c} $,
defined by $ [\underline{c}]\le \underline{c}  <[\underline{c}]+1$ and set
$N=\max \{ N_1, [\underline{c}] \}.$

\begin{itemize}

\item[(H7)]   For any  multi index $\mu \in \mathbb{Z}^{n+1}$
with $|\mu|:=|\mu_0|+...+|\mu_k|\leq 2N +1$ and any $m\in \Z$    with $|m|\le N $,
    we have $\mu \cdot
\lambda +m\neq \underline{c}$.

\item[(H8)] If $0<\lambda _{j_1}<...<\lambda _{j_k}$ are $k$ distinct
  $\lambda$'s,   $\mu\in \Z^k$ satisfies
  $|\mu| \leq 4N +2$   and $m\in \Z $ satisfies  $|m|\le  2N$, then we have
$$
\mu _1\lambda _{j_1}+\dots +\mu _k\lambda _{j_k}+m=0 \iff \mu=0\  \text{ and } m=0.
$$

\item [(H9)]
The Fermi golden rule  Hypothesis (H9')   in subsection
\ref{subsec:FGR}, see \eqref{eq:FGR}, holds.

\item [(H10)] We have $\gamma _1\neq 0$ in \eqref{eq:nonlinearity}.

\end{itemize}

\begin{theorem}\label{theorem-1.1}
    Let $u(t,x)$ be a solution to
\eqref{NLS1}. Assume (H1)--(H10). Then, there exist an $\epsilon_0>0$
and a $C>0$ such that if $  \|u_0  \|_{H^1}<\epsilon  $ with
$\epsilon \in (0, \epsilon_0)$,  there exist
  $h _\pm \in
H^1$ with $\| h_\pm\| _{H^1} \le C \|u_0  \|_{H^1} $ such that

\begin{equation}\label{scattering}
\lim_{t\to  \pm\infty}\|u(t,\cdot) -e^{\im t\Delta }h _\pm\|_{H^1}=0
.
\end{equation}
It is possible to write $u(t,x)=A(t,x)+\widetilde{u}(t,x)$ with
$|A(t,x)|\le C_N(t) \langle x \rangle ^{-N}$ for any $N$, with $\lim
_{|t|\to \infty }C_N(t)=0$ and such that for any pair $(r,p)$ which
is admissible, by which we mean that
\begin{equation}\label{admissiblepair}  2/r+3/p= 3/2\,
 , \quad 6\ge p\ge 2\, , \quad
r\ge 2,
\end{equation}
we have
\begin{equation}\label{Strichartz} \|  \widetilde{u} \|
_{L^r_t( \mathbb{R},W^{1,p}_x)}\le
 C\|  u_0  \| _{H^1 }.
\end{equation}

\end{theorem}

\begin{remark}
\label{rem:h10} When $\gamma _1=0$ equation \eqref{NLS1} admits
standing waves of arbitrarily small energy, by simple bifurcation
theory. So (H10), that is
 $\gamma _1\neq 0$,   is an essential hypothesis.
\end{remark}

\begin{remark}
\label{rem:ckp} Theorem  \ref{theorem-1.1} is a generalization of
the main result of \cite{ckp} which focuses on the special case
$n=0$ and $\underline{c}<1$. In the special case treated in
\cite{ckp}, our proof is particularly simple (although this is here
obscured by our emphasis  on easing the restrictions on $\sigma
(\mathcal{H})$ of \cite{ckp}). Notice that we do not obtain
analogues of the decay formulas (1.12)--(1.14) \cite{ckp} because
our initial data are not required to satisfy $\int _{\R ^3}
|x|^\sigma |u_0(x)|^2dx \ll 1$ for a $\sigma >5$, like in
\cite{ckp}: if one merely asks $\| u_0\| _{H^1}\ll 1$, as we do,
the decay formulas (1.12)--(1.14) \cite{ckp}  are not true.
\end{remark}
\begin{remark}
\label{rem:ckp1} At the beginning of section  7 \cite{ckp}  is
mentioned, without details, the possibility of proving the main
result of \cite{ckp} in the case $n=0$ and $\underline{c}>1$. In
fact this case is substantially harder, even more    if also $n\ge 1$. Treating these cases is
what we accomplish here.
\end{remark}
\begin{remark}
\label{rem:ckp2} We treat a larger class of solutions than
\cite{ckp} and we are able to draw stronger conclusions.   For
instance, in \cite{ckp} the only energy bound proved is of the form
$ \| u(t)\|_{H^1}\lesssim \epsilon \log\langle   \epsilon^4
t\rangle $: here we prove $ \| u(t)\|_{H^1}\lesssim \epsilon
  $.
\end{remark}

\begin{remark}
\label{rem:nonlinearity}  We choose the nonlinearity $|u|^2u$ to
simplify exposition. Indeed in this case the energy $E(t,u)$, see
\eqref{eq:energyfunctional}, is smooth in $(t,u)\in \R \times
H^1(\R ^3)$, and Theorem \ref{theor:normal form} below is easier to
prove. But, with   more effort and with essentially the same
argument, we could have considered a nonlinearity $\beta (|u|^2)u$,
with $\beta (0)=0$, $\beta\in C^\infty(\R,\R)$ and s.t. there
exists a $p\in(1,5)$ s.t. for every $k\ge 0$ there is a fixed $C_k$
with $ \left| \frac{d^k}{dv^k}\beta(v^2)\right|\le C_k |v|^{p-k-1}
\quad\text{if $|v|\ge 1$},$  see \cite{bambusicuccagna,Cu1}.
\end{remark}

\begin{remark}
\label{rem:dimension} We choose space dimension 3 only for
definiteness. It is possible to prove a similar theorem for any
spatial dimension.  In low   dimensions, 1 and 2, there are no
endpoint Strichartz inequalities,  but  \cite{M1,M2} gives us good
  surrogates. Notice that the proofs in \cite{M1,M2} can be substantially simplified, following the ideas from Lemma 3.2 to Lemma 3.6 in \cite{cuccagnatarulli}.
The nonlinearity $|u|^2u$ can be
treated in space dimension 2.   For space
dimension 1, given our need of Strichartz estimates to close
nonlinear estimates, it is necessary to work with $\beta (|u|^2)u$,
$\beta\in C^\infty(\R,\R)$ with $\beta (0)=\beta '(0)=0$ . Hence, and
this is an important technical constraint, the nonlinearity is 0 at
least at fifth order    at $u=0$.
 For 1 D the
nonlinearity $|u|^2u$ is difficult, being {\it long range}, so it
remains an open problem. Notice that for space dimension 1, the
energy $E(t,u)$ is smooth in $(t,u)\in \R \times H^1(\R ^3)$  for
all $\beta (|u|^2)u$ with $\beta\in C^\infty(\R,\R)$, so Theorem
\ref{theor:normal form} below can be easily proved with $\beta
(|u|^2)u$. Summing up,  the second  of the open problems stated at
the end of p.115 \cite{ckp} is mostly solvable, but is unsolved  in
the case of the nonlinearity $|u|^2u$ in  1 D.
\end{remark}

\begin{remark}
\label{rem:genericity} Hypothesis (H9) probably holds for generic
  $ V $. Under the hypotheses on $\sigma (\mathcal{H})$ of
\cite{ckp}, hypothesis (H9) is the same of formula (1.9) in
\cite{ckp}. It is easy to show that (1.9) in \cite{ckp}  holds for
generic $V$, if additionally $|V(x)|\le C e^{-a|x|}$ for $a>0$, see
the proof of Proposition 2.2 and Remark A.1 in
\cite{bambusicuccagna}. Notice that Proposition 2.2
\cite{bambusicuccagna} proves that the   analogue of (H9)  in
\cite{bambusicuccagna} is true for generic $ \beta (|u|^2) $ with
fixed $V$ (with $V$ exponentially decreasing and with simple
eigenvalues): probably an analogous  proof yields  (H9)    for
generic pairs $( \beta (|u|^2) , V)$.
\end{remark}
\begin{remark}
\label{rem:gamma} The function $\gamma (t)$ in
\eqref{eq:nonlinearity} is particularly simple. This simplifies the
exposition. But  similar arguments  work if $\gamma (t)$ is a
higher degree trigonometric polynomial, or if,  for $P(x_1,
y_1,...,x_A, y_A)$ a
 real valued nonconstant polynomial,
$ \gamma(t) = P ( \cos (\omega _1t), \sin (\omega _1t), ..., \cos
(\omega _At), \sin (\omega _At)   ) $, adding appropriate non
resonance hypotheses on these frequencies, the eigenvalues of
$\mathcal{H}$ and $\underline{c}$.
\end{remark}

We recall that  \cite{ckp} shows that (under very restrictive
hypotheses) nonlinear coupling of  continuous and discrete modes
is responsible of
  leaking of energy from discrete modes into radiation,
where linear dispersion occurs.  This is analogous in linear theory
to the Stark effect, see  \cite{Y}, and to effects of disturbances
on ground states, see \cite{KW} and references therein. Nonlinear
coupling       of  continuous and discrete modes is exploited  in
\cite{SW3,bambusicuccagna} for a proof of scattering of small
energy solutions of the nonlinear Klein Gordon equation (NLKG) with
discrete modes. The same idea is exploited   in a substantial
number of papers dealing with asymptotic stability of ground states
of the nonlinear Schr\"odinger equation (NLS), see \cite{Cu1} and
therein for more references.   From the beginning,  at least 15
years ago \cite{sigal,BP2,SW3}, it was clear how coupling should
act. See also the improvements on \cite{BP2,SW3} contained in
\cite{cuccagnamizumachi,zhousigal}. In particular, attention was
drawn to the sign of specific coefficients of the discrete  mode
equations. This sign is responsible for friction on the discrete
modes. Since the system is hamiltonian,  the energy is conserved
and simply is moving from discrete to continuous modes. Except for
special cases though, it  was not  clear     how to prove the sign.
Emphasis was rightly attached to   the necessity  of spectral
resonance between eigenvalues and continuous spectrum of the
linearization. Indeed when this resonance is absent,
 like for the examples of discrete NLS in \cite{Cu2},
  the discrete modes persist. But for continuous NLS there is always
   spectral
resonance. What was not well appreciated was  the crucial role of
the hamiltonian structure. Attempts to prove friction without
exploiting the hamiltonian structure were extremely complex, see
\cite{Gz}, or somewhat indirect  and unsatisfying, see
\cite{cuccagnamizumachi}, and    confined to the case of just 1
discrete mode. In the case of multiple discrete modes, friction was
proved only in special cases with the eigenvalues   close to the
continuous spectrum, \cite{T,Cu3,zhouweinstein1}. The first
reference which recognizes the relevance of the hamiltonian
structure seems to be \cite{Cu3}. In \cite{bambusicuccagna,Cu1} we
were finally able to exploit the intuition of \cite{Cu3} and to
exploit the hamiltonian structure to prove the friction on the
discrete modes under very general   hypotheses on the spectrum.
  By applying a conceptually simple
form of the Birkhoff normal form argument, \cite{bambusicuccagna}
is able to extend the result in \cite{SW3} by dropping   the
spectral assumptions in \cite{SW3}. Along the same lines,
\cite{Cu1}   proves asymptotic stability of ground states of the
NLS. The situation in   \cite{Cu1} is harder, since rather that an
equilibrium point there is an invariant manifold.
  In the present paper we return to
the easier setting of \cite{bambusicuccagna} where the issue is to
prove the asymptotic stability of the equilibrium point 0. So, we
improve \cite{ckp} in the same way  \cite{bambusicuccagna} improves
\cite{SW3}. As in \cite{bambusicuccagna}, the Birkhoff normal form
helps us to prepare the system. We then conclude the proof  in
Section \ref{sec:dispersion}   with somewhat standard arguments,
derived most directly   from
\cite{cuccagnamizumachi,bambusicuccagna,Cu1}, but  which  are a
simplification and generalization of arguments already in
\cite{BP2,SW3}. Notice that these arguments are simpler that
\cite{ckp}. For example, there is no need of hierarchies of Banach
spaces like in \cite{ckp}. Finally, for   the physical relevance,
some interesting open problems and more context and references,  we
refer to \cite{ckp}.

We end the introduction with some notation. Given two functions
$f,g:\mathbb{R}^3\to \mathbb{C}$ we set $\langle f,g\rangle = \int
_{\mathbb{R}^3}f(x)  g(x)  dx$.   For any $k,s\in \mathbb{R}$
 we set
\[ H^{ k,s}(\mathbb{R}^3)=\{ f:\mathbb{R}^3\to \C \text{ s.t.}
 \| f\| _{H^{s,k}}:=\| \langle x \rangle ^s  (-\Delta +1)^{k} f
  \| _{L^2
 }<\infty \}.\]
 We set $  \mathcal{S} (\mathbb{R}^3 )  =\cap _{s,k=0}^{\infty}H^{
k,s}(\mathbb{R}^3) $.
 We set
 $L^{2,s} =H^{0,s}  $, $L^2=L^{2,0}    $,  $H^k=H^{2,0}    $.
Sometimes, to emphasize that these spaces refer to spatial
variables, we will denote them by $W^{k,p}_x$, $L^{ p}_x$, $H^k_x$,
$H^{ k,s}_x$ and $L^{2,s}_x$. For $I$ an interval and $Y_x$ any of
these spaces, we will consider Banach spaces $L^p_t( I, Y_x)$ with
mixed norm $ \| f\| _{L^p_t( I, Y_x)}:= \| \| f\| _{Y_x} \|
_{L^p_t( I )}.$ Given an operator $A$, we will denote by
$R_A(z)=(A-z)^{-1}$ its resolvent. We set
$\mathbb{N}_0=\mathbb{N}\cup \{0 \}$.

\section{Local well posedness and Hamiltonian structure}\label{sec:loc WP}

The first step is the locally well-posed in $H^1( \mathbb{R}^3)$ of
the initial-value problem (\ref{NLS1}), see Theorem 2 \cite{ckp}:

\begin{theorem}\label{th:local}
For every $u_0\in H^1(  \mathbb{R}^3)$ there exists a unique
solution $u(t)$ of the initial-value problem (\ref{NLS1}) defined on
a maximal interval $t\in [0, T_{max})$ such that
\begin{equation} u\in C^1\left([0,T_{max}),H^{-1}\right) \cap
C\left([0,T_{max}),H^1\right),\label{uspace} \end{equation} where if
$ T_{max}<\infty$ then $\lim_{t\to T_{max}}\|u(t)\|_{H^1}=\infty$.
Moreover, $\|u(t)\|_{2}\equiv \|u(0)\|_{2}$, $\forall t\in
[0,T_{max})$, and $u(t)$ depends continuously on the initial data,
i.e. if $\lim_{n \to \infty} u_0^n = u_0$ in $H^1( \mathbb{R}^3)$
then for any closed interval $I\subset [0,T_{max})$ the solution
$u^n(t)$ of the problem (\ref{NLS1}) with initial data $u^n_0$ is
defined on $I$ for sufficiently large $n$, and $\lim_{n \to \infty}
u^n(t) = u(t)$ in $C(I,H^1)$.
\end{theorem}

  The next step is about the hamiltonian nature of \eqref{NLS1}, which is neglected in
    \cite{ckp} but which  in fact  is crucial. So we spend the rest of the section to discuss
  the hamiltonian set up. We have an energy functional

 \begin{equation} \label{eq:energyfunctional}\begin{aligned}&
 E(t,u)=E_K( u)+E_P(t,u)\\&
E_K( u)= \int _{\R ^3}
  \nabla u \cdot \nabla \overline{u} dx  + \underline{c}\int _{\R ^3}
     u \overline{u} dx + \int _{\R ^3}
  V u   \overline{u} dx \\&
E_P(t,u)=\gamma (t)
 \int _{\R ^3}\frac{|u |^4}{4} dx \end{aligned}
\end{equation}
   Equation
\eqref{NLS1} can be written as
\begin{equation}\label{eq:NLShamiltonian} \im \dot u =
 \partial _{\overline{u}}E(t,u)
  .
\end{equation}
   We consider
  eigenfunctions $\phi _j(x)$   with eigenvalue   $\lambda
_j $: $ \mathcal{H} \phi _j =\lambda _j \phi _j .$ They can be
normalized so that they are real valued and $\langle   \phi  _j ,
{\phi} _\ell  \rangle =\delta _{j\ell }$. The
  $\phi _j( x)$ are smooth  and satisfy for some fixed $a>0$ and all multi indexes $\alpha$
\begin{equation}\label{eq:regdecEigenf}\sup_{ x\in\R^3, j=0,n} e^{a|x|}  |\partial ^\alpha _x\phi_j(x)|
  <\infty .
\end{equation}
We have the $\mathcal{H} $  decomposition
\begin{align}  \label{eq:spectraldecomp} &
L^2(\R^3,\C )= \ker (\mathcal{H})\oplus _{j=1}^n \ker (\mathcal{H}-
\lambda _j ) \oplus L_c^2(\mathcal{H} ).
\end{align}
    Correspondingly we set,
\begin{equation}\label{eq:coordinate}u=  z \cdot \phi   +   f,    \text{ for } z \cdot \phi
=\sum _{j=0}^n z_j \phi _j(x). \end{equation}   In $L^2(\R^3,\C ^2
)$
 we consider the
symplectic form
\begin{equation}
\label{eq:Omega0} \Omega _0((u ,\overline{u}),(v,\overline{v})
):=\langle u ,\overline{v} \rangle - \langle \overline{u} ,v \rangle
\end{equation}
The proof of  the following lemma is elementary:
\begin{lemma}
  \label{lem:Omega _0} Let $P_c$ be the projection in
   $L^2_c(\mathcal{H})$.  Consider the functions $z_j= \langle u,
   \phi _j \rangle $,  $\overline{z}_j= \langle \overline{u},
   \phi _j \rangle $, $f= P_cu$ and $\overline{f}=P_c\overline{u}$.
   Then in terms of these functions we have
\begin{equation}
\label{eq:Omega0coord} \Omega _0= \sum _{j=0}^n dz_j\wedge
d\overline{z}_j+\langle f'  \quad , \overline{f}'  \quad \rangle -
\langle \overline{f}'  \quad ,  {f}'  \quad \rangle ,
\end{equation}
where $f'v=P_cv$ and  $\overline{f}'v=\overline{P_cv}$ and where
$\langle f'  \quad , \overline{f}'  \quad \rangle$ (resp. $ \langle \overline{f}'  \quad ,  {f}'  \quad \rangle$) acts on a pair $(v_1,v_2)$
as  $\langle f'   v_1 , \overline{f}'   v_2 \rangle$
(resp. $ \langle \overline{f}'   v_1 ,  {f}'   v_2 \rangle$).

\end{lemma}

We consider now two additional variables $(t,\tau )\in \R ^2$ and we
set
\begin{equation}
\label{eq:Omega} \Omega =\Omega _0+\im \, dt \wedge d\tau .
\end{equation}
For a function $F$ we call   hamiltonian vector field $X_F$ with respect
to $\Omega $ the field defined by $\Omega ( X_F,Y)=-\im \, dF (Y)$. For any
vector $\Sigma $ we set
\begin{equation} \label{eq:ExpVect}  \begin{aligned} &
\Sigma =\Sigma  _{t}\frac{\partial}{\partial t}+\Sigma
_{\tau}\frac{\partial}{\partial \tau} +\sum \Sigma
_{j}\frac{\partial}{\partial z_j}+\sum \Sigma
_{\overline{j}}\frac{\partial}{\partial \overline{z}_j} + \Sigma _{f
} +\Sigma _{\overline{f} }
\end{aligned}\end{equation}
for
\begin{equation} \label{eq:Y1}  \begin{aligned} & \Sigma _{t}
=dt (\Sigma )\, , \quad  \Sigma _{\tau} =d\tau (\Sigma )\, , \quad
\Sigma_{j} =dz_j (\Sigma)\\& \quad \Sigma _{\overline{j}}
=d\overline{z}_j (\Sigma ) \, , \quad  \Sigma_{f } =f' \Sigma \, ,
\quad  \Sigma_{\overline{f} } =\overline{f}' \Sigma.
\end{aligned}\end{equation}
A  differential 1-form $\alpha $ decomposes as
\begin{equation} \label{eq:gamma}  \begin{aligned} &
\alpha=\alpha ^{t}d t +\alpha ^{\tau}d \tau +\sum \alpha ^{j}d z_j
+\sum \alpha ^{\overline{j}}d\overline{z}_j +\langle \alpha ^{f },
f'\quad \rangle +\langle \alpha ^{\overline{f} }, \overline{f}'\quad
\rangle   ,
\end{aligned}\end{equation}
where $\langle \alpha ^{f },
f'\quad \rangle  $ (resp. $\langle \alpha ^{\overline{f} }, \overline{f}'\quad
\rangle $) acts on a vector $v$ as $\langle \alpha ^{f },
f'v \rangle  $ (resp. as $\langle \alpha ^{\overline{f} }, \overline{f}'v
\rangle $).
Then
\begin{equation}
\label{eq:hamiltonian field} \begin{aligned} & (X_F)_t =-
\frac{\partial
  F}{\partial \tau} \, , \, (X_F)_\tau =
\frac{\partial
  F}{\partial t} \, , \, (X_F)_j=-\im\frac{\partial
  F}{\partial \overline{z}_j}\\&
  (X_F)_{\overline{j}}=\im \frac{\partial F}{\partial z_j} \, , \,
  (X_F)_{f}=-\im \nabla_{  \overline{f}}F \, , \,
  (X_F)_{\overline{f}}= \im \nabla_{   {f}}F.\end{aligned}
\end{equation}
We call   Poisson bracket of two functions with respect to $\Omega $
the function
\begin{equation}
\label{eq:PoissBrack} \begin{aligned} &\left\{F,G\right\}:= dF
(X_G)=\frac{\partial F}{\partial
  \tau }\frac{\partial G}{\partial t}-\frac{\partial F}{\partial
  t}\frac{\partial G}{\partial  \tau}+\\& +
  \im \sum_{j=1}^{n}\left(\frac{\partial F}{\partial
  \overline{z}_j}\frac{\partial G}{\partial z_j}-\frac{\partial F}{\partial
  z_j}\frac{\partial G}{\partial  \overline{z}_j} \right)+  \im \left\langle
  \nabla_{  \overline{f}}F,\nabla_fG \right\rangle- \im \left\langle
  \nabla_{  \overline{f}}G,\nabla_fF \right\rangle .\end{aligned}
\end{equation}
We consider the new hamiltonian

\begin{equation}\label{eq:H} H:= H_F+     E_P(t, u)\, , \quad H_F:=
E_K-\tau .
\end{equation}
 Then the corresponding hamiltonian system is
\begin{equation} \label{eq:SystemH} \begin{aligned} &
   \frac{\partial}{\partial s} z_j  =- \im  \frac{\partial H}{\partial   \overline{z}_j  }
\, , \quad  \frac{\partial}{\partial s} f=  - \im \nabla
_{\overline{f}}  H\\& \frac{d}{d s} t=-  \frac{\partial}{\partial
\tau}H =1\, , \quad  \frac{d}{d s} \tau =   \frac{\partial}{\partial
t}H .
\end{aligned}
\end{equation}
It is easy to see that \eqref{eq:SystemH} and
\eqref{eq:NLShamiltonian} are equivalent.   Plugging the decomposition \eqref{eq:coordinate} in the energy \eqref{eq:energyfunctional}, is is easy to see the following equality:
\begin{equation}\label{eq:HF}   H_F =\sum _{j=1}^n\lambda _j
z_j \overline{z}_j  +    \langle \mathcal{H}f, \overline{f}\rangle
-\tau .
\end{equation}

Once the natural coordinates \eqref{eq:coordinate} are introduced,
\cite{ckp} starts a normal form argument to simplify the system.
Here we do the same, but we want to preserve the hamiltonian
structure. So  we need canonical transformations, which is the theme
of   section \ref{section:normal forms}.

\section{Canonical transformations}
\label{section:normal forms}

\subsection{Lie transform}
\label{subsec:LieTransf}   For $ {m}_0 \in \mathbb{N}_0$, $M_0\in
\mathbb{N}$ we consider functions

\begin{equation}
\label{chi.1}  \begin{aligned}   \chi  =  \sum  _{ \ell = - {m}_0
}^{m_0}  e^{\im \ell t}  [ & \sum_{\substack{
  |\mu |=|\nu |=M_0+1}} a_{\ell\mu \nu }z^{\mu}
\overline{z}^{\nu} + \sum_{\substack{|\mu   |=M_0 \\ |\mu |=|\nu
|-1}}  z^{\mu} \overline{z}^{\nu}
 \langle    \Phi _{\ell\mu   \nu
}
  , f \rangle \\& +  \sum_{\substack{| \nu |=M_0 \\ |\mu |=|\nu |+1}}  z^{\mu} \overline{z}^{\nu}
 \langle    \Psi _{\ell\mu   \nu
}
  , \overline{f} \rangle ] .\end{aligned}\end{equation}
  We assume
\begin{equation} \label{eq:symmetries}\begin{aligned}&
 \overline{a}_{\ell\mu \nu
}= a_{-\ell \nu \mu }\, , \quad \overline{\Phi}_{\ell \mu \nu }(x)=\Psi
_{-\ell \nu \mu }(x).
\end{aligned}
\end{equation}
We assume    $\Phi_{m\mu   \nu }(x)\in  \mathcal{S} (\mathbb{R}^3 )
$.  Denote by $\mathfrak{F}  ^s$ the flow of the Hamiltonian vector
field $X_{\chi}$. The {\it Lie transform} $\mathfrak{F} =
 \mathfrak{F}^s\big|_{s=1}$ is defined in a sufficiently small neighborhood
 of the origin and is a canonical transformation. Then
 we have, with $(z',f')=\phi (
 t ,z,f)$,
 \begin{equation} \label{eq:Lie}(\tau ',
 t',z',f'):=\mathfrak{F}(\tau ,t,z,f)=  (\tau +\psi
 (t,z,f)   ,t+1,\phi ( t,z,f))  .\end{equation}

\begin{lemma}\label{lie_trans}
Consider the $\chi$ in \eqref{chi.1} and its Lie transform $\mathfrak{F} $.
  Fix any $n_0\in \mathbb{N}\cup \{ 0 \}$.
 Then there are $ \mathcal{G}(t ,z,f )$ and    $\Gamma _j (t ,z,f )$
 such that for any pair $(-K',-S')$
we have  the following properties.
\begin{itemize}
\item[(1)]    $\Gamma _j \in
C^\infty ( \U ^{-K',-S'}, \mathbb{R}) $, with $\U ^{-K',-S'}=\widetilde{\U} ^{-K',-S'}\times \R $,
with $\widetilde{\U} ^{-K',-S'}
\subset
\mathbb{C}^{n+1} \times H^{-K',-S'}_c   $ an
appropriately small neighborhood of the origin.
\item[(2)] $  \mathcal{G}\in C^\infty
(\U ^{-K',-S'} ,  H^{K,S}_c  )  $  for any $K,S$.
\item[(3)] The transformation $\phi (t ,z,f)$ is such that we have

\begin{equation} \label{eq:LieTransf}\begin{aligned}&    z' _j =
  z _j+U_j+\Gamma _j  (t ,z,f ) ,\\
  & f' =  f +V + \mathcal{G}   (t ,z,f )
  \end{aligned}\end{equation}
where
\begin{equation} \label{eq:LieTransf1}\begin{aligned}& U_j=
\sum _{l=1}^{n_0 } \frac{1}{l!}  \lie _\chi ^l(z_j)   ,  \quad  V=
\sum _{l=1}^{n_0  } \frac{1}{l!}  \lie _\chi ^l(f),  \quad   \lie
_\chi ^l(g):=\{ \{ .. \{ g, \underbrace{\chi \} ..  \chi } _{l}\}
\\& \{f, \chi \} :=-\im \nabla _{\overline{f}}\chi
 =
-\im \sum  _{\ell = -{m}_0 }^{m_0} e^{\im \ell t} \sum_{\substack{| \nu |=M_0
\\ |\mu |=|\nu |+1}}  z^{\mu} \overline{z}^{\nu}
     \Psi _{\ell\mu   \nu} (x)\, , \\&  \lie _\chi ^l(f): = -
     \im \sum  _{\ell = -{\ell}_0 }^{m_0} e^{\im \ell t}   \sum_{\substack{| \nu |=M_0
\\ |\mu |=|\nu |+1}}\lie _\chi  ^{l-1 }(
z^{\mu} \overline{z}^{\nu}) \Psi _{\ell\mu   \nu} (x).
 \end{aligned}\end{equation}

\item[(4)]  $\lie _\chi ^l(z_j)$ is a  homogeneous polynomial    of degree  $2l
 M_0        +1$ in $(z,\overline{z})$, $\langle    \Phi _{\ell\mu \nu
}
  , f \rangle$ and $\langle    \Psi _{\ell\mu   \nu
}
  , \overline{f} \rangle $. Same statement is true for
    $\lie _\chi ^l(f)$.

\item[(5)] There exists a constant $C=C(K,S,K',S')$ such that near 0
\begin{equation} \label{eq:polEst}\begin{aligned}& |\Gamma _j|  +
\| \mathcal{G}\| _{H^{K,S}}   \le C  (|z| +\norma{f} _{H^{-K',-S'}}
) ^{   2( n_0+1) M_0  +1}  .\end{aligned}
\end{equation}

\end{itemize}
\end{lemma}
\proof Recall that for any function $\psi$, we have $\frac{d}{ds}
(\psi \circ \mathfrak{F}^s)=\{ \psi, \chi \} \circ \mathfrak{F}^s.$
Then
   \eqref{eq:LieTransf} follow
by Taylor expansion with reminders

\begin{equation}\label{eq:resti}\begin{aligned} & \Gamma
_j:=
 \int _0^1  \frac{(1-s )^{{n_0  }}} {{ n_0  }!} \{
\{.. \{ z_j, \underbrace{\chi \} .. \chi } _{n_0 +1 } \} \circ
\mathfrak{F} _s ds  \\& \mathcal{G}:= \int _0^1 \frac{(1-s )^{{n_0
}}} {{ n_0  }!} \{ \{...\{ f, \underbrace{\chi \} ...\chi } _{n_0 +1
} \} \circ \mathfrak{F} _s ds.
\end{aligned}
\end{equation}
\eqref{eq:polEst} follows closing up inequalities with a standard
argument, see proof of Lemma 4.3 \cite{bambusicuccagna}. Claim (4) is true for $l=1$, see for example $\lie _\chi  (f) $
in \eqref{eq:LieTransf1}. The general case follows by induction. \qed

\begin{lemma}\label{lem:lie_trans0}    Consider\begin{equation} \label{eq: LieTransf g0} g=e^{\im t m}
z^{\mu}\overline{z}^{\nu} \langle    \Phi
  , f \rangle  ^{\alpha}\langle    \Psi
  , \overline{f} \rangle
  ^{\beta} \langle f^{a}\overline{f} ^{b}, \Psi _{ab}\rangle ,  \end{equation} where for a multiindex $\alpha =( \alpha _1,..., \alpha _M)$ we have set
  $\langle    \Phi
  , f \rangle  ^{\alpha} =\prod _{j=1}^{M}\langle    \Phi _j
  , f \rangle  ^{\alpha _j}$ for     $\Phi _{j} \in \mathcal{S}(\R ^3, \C) $
  and where $ \langle    \Psi
  , \overline{f} \rangle
  ^{\beta} $ has a similar meaning, i.e. $ \langle    \Psi
  , \overline{f} \rangle
  ^{\beta} =\prod _{j=1}^{\widetilde{M}}\langle    \Psi _j
  , \overline{f} \rangle  ^{\beta _j}$ for     $\Psi _{j} \in \mathcal{S}(\R ^3, \C) $
  for a multiindex $\beta =( \beta _1,..., \beta  _{\widetilde{M}})$.
  Assume:

\begin{itemize}
\item[(i)]
  $a+b\le 4$;

  \item[(ii)]
  $|\mu |+|\alpha |+a=|\nu |+|\beta |+b =L+1;$

  \item[(iii)]
     $\Psi _{ab} \in \mathcal{S}(\R ^3, \C) $ for  $a+b< 4$; if $a+b= 4$  we assume  $\Psi _{ab}= \frac{1}{4}\delta _{ab}$ and $L=2$;

      \item[(iv)] $|m|\le L.$

\end{itemize}
      Consider $\chi$ as in \eqref{chi.1}. Then   $\lie _\chi  (g)$ is a finite sum   of  terms of the form
  \begin{equation} \label{eq: LieTransf g1}  e^{\im tm'
  }
z^{\mu '}\overline{z}^{\nu '} \langle    \Phi '
  , f \rangle  ^{\alpha '}\langle    \Psi '
  , \overline{f} \rangle
  ^{\beta '} \langle f^{a'}\overline{f} ^{b'}, \Psi ' _{a'b'}\rangle  \end{equation}
with:
\begin{itemize}
\item[(1)]  $a'+b'< 4$ and   $\Psi '_{a'b'} \in \mathcal{S}(\R ^3, \C) $;

\item[(2)] $ |\mu '|+|\alpha '|+a'=|\nu '|+|\beta '|+b'=L'+1  $, with  $L'=L+ M_0 $;

    \item[(3)] $|m'|\le m_0+|m|.$

\end{itemize}
   In \eqref{eq: LieTransf g1}  the factors
$\langle    \Phi '
  , f \rangle  ^{\alpha '}$ and $\langle    \Psi '
  , \overline{f} \rangle
  ^{\beta '} $ are defined like the corresponding ones in \eqref{eq: LieTransf g0} and  involve functions  $\Phi _{j}',\Psi _{j}' \in \mathcal{S}(\R ^3, \C) $.
\end{lemma}
\proof   $\lie _\chi  (g)$ is obtained applying the Leibnitz rule to the rhs\eqref{eq: LieTransf g0}.  Trivially, if $a+b<4$ then also $a'+b'< 4$. On the other hand, for   $a+b=4$ we have $g=\frac{1}{4} \| f \| ^4_4 $ by (iii).
Then $\lie _\chi  (g)$ is a linear combination of terms $\langle  \lie _\chi  (f)  f \overline{f}^2, 1 \rangle $ and $\langle  \lie _\chi  (\overline{f})  \overline{f}  {f}^2, 1 \rangle $. If we apply the formula for $\lie _\chi  (f) $ in \eqref{eq:LieTransf1} and by the fact that $\lie _\chi  (\overline{f}) =\overline{\lie _\chi  (f)} $, consequence of \eqref{eq:symmetries}, we obtain claim (1). (3) follows by the fact that in \eqref{eq: LieTransf g1}  we have
$m'=\ell +m$ with $|\ell |\le m_0$ and $m$ the same of \eqref{eq: LieTransf g0}. Finally, (2) is an elementary consequence
of the   Leibnitz rule and the definition of $\chi $ in \eqref{chi.1}.
For example, one of the terms in the expansion of  $\lie _\chi  (g)$ is
\begin{equation}  ae^{\im t m}
z^{\mu}\overline{z}^{\nu} \langle    \Phi
  , f \rangle  ^{\alpha}\langle    \Psi
  , \overline{f} \rangle
  ^{\beta} \langle \lie _\chi  (f)f^{a-1}\overline{f} ^{b}, \Psi _{ab}\rangle .  \nonumber \end{equation}
If we substitute the formula for $\lie _\chi  (f) $ in \eqref{eq:LieTransf1}
we get  a linear combination of  terms

\begin{equation}   e^{\im t (m+\ell )}
z^{ {\mu + \mu '' } }\overline{z}^{\nu+\nu '' } \langle    \Phi
  , f \rangle  ^{\alpha}\langle    \Psi
  , \overline{f} \rangle
  ^{\beta} \langle  f^{a-1}\overline{f} ^{b},\Psi _{\ell \mu '' \nu '' } \Psi _{ab}\rangle ,  \nonumber \end{equation}
where $|\nu '' |=M_0$,
 $ |\mu '' |=M_0+1$ and $\Psi _{\ell \mu '' \nu '' } \Psi _{ab} \in \mathcal{S}(\R ^3, \C).$ This is of form \eqref{eq: LieTransf g1} with

\begin{equation}
 \begin{aligned} & m'=m+\ell \, , \quad \mu '= \mu +\mu '' \, , \quad \nu '= \nu +\nu ''  \, , \quad  \alpha '=\alpha  \, , \\&   \beta  '=\beta \, , \quad a'=a-1\, , \quad  b'=b \, , \\&  |\mu '|+|\alpha '|+a'= |\mu  |+|\alpha |+a-1+|\mu ''| =
     L+M_0+1, \\&  |\nu '|+|\beta '|+b'= |\nu  |+|\beta |+b +|\nu ''| =
     L+M_0+1.
\end{aligned}\nonumber
\end{equation}

\qed

\begin{lemma}\label{lem:lie_trans1}
  Consider $g$ like in \eqref{eq: LieTransf g0}, with the above conventions
  and with assumptions (i)--(iv) of Lemma \ref{lem:lie_trans1}.

   \begin{itemize}
\item[(a)] We have an expansion
   \begin{equation} \label{eq:LieTransf g}\begin{aligned}&
g\circ \mathfrak{F}_1=g
   +\mathcal{V} +  \mathbf{{G}}, \quad \mathcal{V}
   =\sum _{l=1}^{n_0 } \frac{1}{l!}  \lie _\chi ^l(g),\\&
    |\mathbf{{G}}|\le C (|z| +\norma{f} _{H^{-K',-S'}}
) ^{2 M_0   ( n_0+1)+2L-a-b) } (|z| +\norma{f} _{H^{1}} ) ^{a+b }.
  \end{aligned}\end{equation}
 \item[(b)] Each  $\lie _\chi ^l(g)$ is a finite sum   of  terms of the form \eqref{eq: LieTransf g1} with    $|m'|\le l m_0+|m|$, where  $L'+1=|\mu '|+|\alpha '|+a'=|\nu '|+|\beta '|+b'  $
satisfies $L'=L+l   M_0 $.   We have $a'+b'< 4$ and   $\Psi '_{a'b'} \in \mathcal{S}(\R ^3, \C) $. \end{itemize}
\end{lemma}
 \proof  (b) is obtained iterating the result in Lemma \ref{lem:lie_trans1}.
 We turn to the proof of (a).
By Lemma \ref{eq:LieTransf}, $f\circ \mathfrak{F}_s - f\in
C^{\infty}(\U ^{-K',-S'}, H_c^{K,S})$. Hence   $g\circ
\mathfrak{F}_s$  is smooth. We can apply $\frac{d}{ds} (g \circ
\mathfrak{F}^s)=\{ g, \chi \} \circ \mathfrak{F}^s $ to a  Taylor
expansion with reminder

\begin{equation}\label{eq:resto}\begin{aligned} & \mathbf{{G}}:=
 \int _0^1  \frac{(1-s )^{{n_0  }}} {{ n_0  }!}
 \lie _\chi ^{n_0 +1 } ( g)   \circ
\mathfrak{F} _s ds   .
\end{aligned}
\end{equation}
$\lie _\chi ^{n_0 +1 } ( g)$ is a linear combination for $k+\ell
=n_0 +1  $   of
\begin{equation} \label{eq:rest0}\begin{aligned} & \lie _\chi ^{\ell }(
z^{\mu}\overline{z}^{\nu} \langle    \Phi
  , f \rangle  ^{\alpha}\langle    \Psi
  , \overline{f} \rangle
  ^{\beta} )  \langle \lie _\chi ^{k }(f^{a}
  \overline{f} ^{b}), \Psi _{ab}\rangle .
\end{aligned}
\end{equation}
By Lemma \ref{lie_trans} we have

\begin{equation} \label{eq:rest1} |\lie _\chi ^{\ell }(
z^{\mu}\overline{z}^{\nu} \langle \Phi
  , f \rangle  ^{\alpha}\langle    \Psi
  , \overline{f} \rangle
  ^{\beta} ) | \le C (|z| +\norma{f} _{H^{-K',-S'}}
) ^{2 M_0   ( \ell+1)+2L-a-b}  \end{equation}
and
\begin{equation} \label{eq:rest2}\begin{aligned} &
\langle \lie _\chi ^{k }(f^{a}
  \overline{f} ^{b}), \Psi _{ab}\rangle =\sum _{\sum k_j+\sum k_j'= k} c _{k_1...k_ak'_1...k_b}
 \langle  \prod _{j=1}^{a}\lie _\chi ^{k _j }(f  )
 \prod _{j=1}^{b}\lie _\chi ^{k _j' }( \overline{f} ) , \Psi _{ab}\rangle .
\end{aligned}
\end{equation}
It is easy to conclude that the largest terms in \eqref{eq:rest0} are those with $k=0$
and $\ell =n_0+1$. This yields \eqref{eq:LieTransf g}. When $g=\langle
f^{a}\overline{f} ^{b}, \Psi _{ab}\rangle$, the worst terms are the ones of
the form
\begin{equation} \label{eq:rest3}\begin{aligned} &
|\langle \lie _\chi ^{n_0 +1}(f ) f^{a-1} \overline{f} ^{b}
   , \Psi _{ab}\rangle  |\le C \sum_{\substack{|\mu
  |=M_0
\\ |\mu |=|\nu |+1}} |\lie _\chi  ^{n_0  }(
z^{\mu} \overline{z}^{\nu})|  \| f\| _{H^1}^{a+b-1} \le \\& \le C\|
f\| _{H^1}^{a+b-1}(|z| +\norma{f} _{H^{-K',-S'}} ) ^{ 2 M_0
n_0+2M_0+1}  \\& \le C(|z| +\|
f\| _{H^1})^{a+b }(|z| +\norma{f} _{H^{-K',-S'}} ) ^{2 M_0
(n_0+1)} .
\end{aligned}
\end{equation}

 \qed

\subsection{Normal forms}
  We set $\lambda = (0, \lambda _1, ..., \lambda _n)$.
\begin{definition}
\label{def:normal form} A function $Z(t,z,
\overline{z},f,\overline{f})$ is in normal form if it is of the form
\begin{equation}
\label{e.12} Z(t,z, \overline{z},f,\overline{f})=Z_0( z,
\overline{z} )+Z_1(t,z, \overline{z},f,\overline{f})
\end{equation}
where we assume properties (N0)--(N2)   listed now.
\begin{itemize}
\item[(N0)] $Z_0( z,
\overline{z} )$ is a finite sum

\begin{eqnarray}
\label{eq:Z0} Z_0( z, \overline{z} )= \sum _{\mu \nu} a_{\mu   \nu}
 z^\mu \overline{z}^\nu \text{ with $a_{\mu   \nu}\in \C$ and with}\\
\label{eq:Z01}  \lambda  \cdot(\mu-\nu)=0 \text{ and }|\mu| =|\nu |
.
\end{eqnarray}

\item[(N1)] $Z_1(t,z, \overline{z},f,\overline{f})$ is a finite sum
\begin{equation}
\label{eq:Z1} \sum _{m ,\mu ,\nu}e^{\im tm} z^\mu \overline{z}
^\nu\langle \Phi _{m  \mu  \nu}, f\rangle + \sum _{m ',\mu ' ,\nu
'}e^{\im tm'} z^{\mu'}\overline{z}^{\nu'}\langle \Psi _{m ' \mu '
\nu '}, \overline{f}\rangle
\end{equation}
with  $\Phi _{m   \mu  \nu}(x) ,\Psi _{m  ' \mu  ' \nu'} (x)\in
\mathcal{S} (\R ^3, \C )$   and  indexes  satisfying
\begin{eqnarray} &
\label{eq:Z11}  \lambda \cdot(\mu-\nu)-m  <-\underline{c} \ ,\quad
|\mu | =|\nu |-1  \ ,\quad |m|\le |\mu |
\\& \label{eq:Z12}  \lambda \cdot (\mu'-\nu')-m'>\underline{c} \ ,\quad |\mu '| =|\nu '
 |+1  \ ,\quad |m|\le |\nu '|  .\end{eqnarray}
\item[(N2)]   $Z_0( z, \overline{z} )$ and $Z_1(t,z, \overline{z},f,\overline{f})$  are real
valued when $\overline{z}$ (resp. $\overline{f}$) is the complex
conjugate of $ {z}$ (resp. $ {f}$), that is
\begin{equation}
\label{eq:nf real} \overline{a}_{\mu   \nu}=a_{\nu   \mu}\, , \quad
\Psi _{m \mu  \nu}(x)= \overline{\Phi} _{-m \nu  \mu}(x).
\end{equation}

\end{itemize}
\end{definition}

We will use  the following table of formulas:

 \begin{equation} \label{PoissBra2} \begin{aligned} & \{
H_F,e^{\im tm} z^{\mu} \overline{z}^{\nu}   \} = \im \left ( \lambda
\cdot (\mu - \nu) -m  \right ) e^{\im tm} z^{\mu} \overline{z}^{\nu}
,\\& \{ H_F, e^{\im tm} z^{\mu} \overline{z}^{\nu} \langle   \Phi
,f\rangle \} = \im e^{\im tm} z^{\mu} \overline{z}^{\nu} \langle
\left ( \mathcal{H} - \lambda \cdot (\nu - \mu) -m   \right )\Phi, f
\rangle ,
\\& \{ H_F, e^{\im tm'} z^{\mu '} \overline{z}^{\nu ' }
\langle   \Psi  ,\overline{f}\rangle \} = - \im e^{\im tm'} z^{\mu'}
\overline{z}^{\nu'} \langle \left ( \mathcal{H} -\lambda \cdot (\mu
'- \nu') +m'    \right )
  \Psi ,\overline{f} \rangle .
\end{aligned}
\end{equation}
We set
\begin{equation} \label{eq:resolvent1} R_{m\mu \nu}=R_{\mathcal{H}}
( \lambda \cdot (\mu - \nu) -m)
\end{equation}
 An immediate consequence of the table \eqref{PoissBra2} is the
 following lemma.

\begin{lemma}
\label{sol.homo} Consider finite sums in $m$, $\mu$ and $\nu$

\begin{equation}
\label{eq.stru.1} K =\sum  e^{\im tm} z^{\mu} \overline{z}^{\nu}
\left ( k_{m \mu\nu}+
 \langle \Phi _{m \mu   \nu
}
  , f \rangle   +
 \langle     \Psi _{m\mu   \nu
}
  , \overline{f} \rangle \right ).
\end{equation}
 Suppose that  $k_{0\mu\nu}=0$  if  $(  \mu ,\nu)$
 satisfies  \eqref{eq:Z01},  $\Phi_{m \mu\nu}=0$ if $(m , \mu ,\nu)$
 satisfies    \eqref{eq:Z11}, $\Psi_{m'\mu'\nu'}=0$ if
 $(m ', \mu ',\nu')$
 satisfies   \eqref{eq:Z12}.   Consider
\begin{equation}
\label{solhomo} \begin{aligned} & \chi = \im  \sum e^{\im tm}
z^{\mu} \overline{z}^{\nu} \big [ \left \langle
    R_{-m\nu \mu}  \Phi_{m \mu   \nu }
  , f \right \rangle +\\&
\frac{  k_{m\mu\nu} }{
  \left ( \lambda \cdot (\mu - \nu) -m
  \right )}    -
    \left \langle
   R_{ m \mu  \nu  }  \Psi_{m  \mu   \nu }
  , \overline{f} \right \rangle \big ]  .\end{aligned}
\end{equation}
Then we have
\begin{equation}
\label{eq:homologicalEq} \left\{ \chi ,H_F\right\}   =K .
\end{equation}
If the coefficients of $K$ satisfy \eqref{eq:nf real}, that is
\begin{equation}
\label{eq:real} \begin{aligned} & \overline{k}_{m \mu\nu}=k_{-m \nu
\mu } \, , \quad \overline{\Phi}_{m \mu \nu }= \Psi_{-m \nu \mu },
\end{aligned} \end{equation}
then also $\chi$ satisfies analogous equalities.
\end{lemma}

Finally, we are able to discuss the main result of section \ref{section:normal forms}.

\begin{theorem}
\label{theor:normal form}   For any integer $N+1\ge r\ge 1$   there
are a neighborhood $ \mathcal{U} $ of the origin in    $H^{1 } (\R
^3, \C)$  and a smooth canonical transformation $\Tr_r: \mathcal{U}
\to H^{1 } (\R ^3, \C)$ s.t.
\begin{equation}
\label{eq:bir1} H^{(r)}:=H\circ \Tr_{r }= H_F+Z^{(r)}+\resto^{(r)}.
\end{equation}
     where:
\begin{itemize}
\item[(i)]$Z^{(r)}$ is in normal form, i.e. satisfies (N0)--(N2),
 with  $2r$--degree  mono-mials;
\item[(ii)]     $\Tr _r = \Tr _{r-1}\circ  \mathfrak{F} _r$,
with $\Tr _{1}$ the identity and $\mathfrak{F} _r$ a transformation
as in  \S \ref{subsec:LieTransf} arising from a polynomial $\chi _r$
as in \eqref{chi.1};
\item[(iii)] we have $\resto^{(r)} = \sum _{d=0}^7\resto^{(r)}_d$
with the following properties:
\begin{itemize}
\item[(iii.0)] we have a finite sum

\begin{equation} \label{eq:resto0}\resto^{(r)}_0(t,z)=
\sum _{\substack{|\mu |=|\nu |  \ge r+1\\   |m|\le |\mu |} } a_{m\mu \nu
}^{(r)}
e^{\im m t}z^\mu \overline{z}^\nu
\end{equation}
 with  $a_{m\mu \nu
}^{(r)}\in \C $ s.t.  $\overline{a_{m\mu \nu
}^{(r)}}= a_{-m\nu \mu
}^{(r)} $;

\item[(iii.1)]  we have a finite sum   $ \resto^{(r)}_1(t,z,f)=$
\begin{equation} \label{eq:resto1} \begin{aligned} &
\sum _{\substack{|\mu
|=|\nu |-1  \ge r\\   |m|\le |\mu |}  }
e^{\im m t} z^\mu \overline{z}^\nu \langle    \Phi
_{m\mu \nu }^{(r)}, f\rangle  + \sum _{\substack{|\nu
|=|\mu |-1  \ge r
\\   |m|\le |\nu |}    } e^{\im m t} z^\mu \overline{z}^\nu \langle \Psi _{m\mu \nu }^{(r)},
\overline{f}\rangle
\end{aligned}
\end{equation}
with  $ \Phi
_{m\mu \nu }^{(r)} \in \mathcal{S} (\R ^3, \C )$ and with  $ \Psi
_{-m\nu \mu }^{(r)} = \overline{\Phi
_{m\mu \nu }^{(r)}}$;

\item[(iii.2--5)] for $d=2,...,5$  we have $\resto^{(r)}_d(t,z,f)= $

\begin{equation}  \begin{aligned} &
\sum _{\substack{a+b=d\\|\mu |+|\alpha |+a=\\  |\nu |+|\beta |+b =:L+1\ge 2\\|m|\le L }} e^{\im m t}z^{\mu}\overline{z}^{\nu} \langle    \Phi ^{(r)}
  , f \rangle  ^{\alpha}\langle    \Psi ^{(r)}
  , \overline{f} \rangle
  ^{\beta} \langle f^{a}\overline{f} ^{b}, \Psi _{abm\mu \nu \alpha\beta}^{(r)}\rangle,\end{aligned} \nonumber
\end{equation} with   functions $ \Phi ^{(r)}_j$, $\Psi ^{(r)} _j$ and $\Psi _{abm\mu \nu \alpha\beta}^{(r)}$ in $\mathcal{S} (\R ^3, \C )$
and with  $\resto^{(r)}_d  $ real valued;

\item[(iii.6)] $ \resto^{(r)}_6(t,z,f)=\gamma (t)
 \int _{\mathbb{R}^3}    | f (x)| ^4 dx/4$;

\item[(iii.7)]  we have $\resto^{(r)}_7\in C^{\infty}(\U ^{-K',-S'}, \R)$ for all $(K',S')$ with    $\U ^{-K',-S'}=\widetilde{\U} ^{-K',-S'}\times \R $,
with $\widetilde{\U} ^{-K',-S'}
\subset
\mathbb{C}^{n+1} \times H^{-K',-S'}_c   $ an
appropriately small neighborhood of the origin, and with
$$ |\resto^{(r)}_7(t,z,f)|\le  C    (|z| +\norma{f} _{H^{-K',-S'}}
) ^{   2N+4}.$$

\end{itemize}
\end{itemize}
\end{theorem}
\proof Case $r=1$  is true with $Z^{(1)}=0$, $\resto ^{(1)}= E_{P}$
and $|m|\le 1$. We proceed by induction. Set

\begin{equation}
\label{pr.3} \begin{aligned}   &\widetilde{K}_{r+1}:=\widetilde{\resto} ^{(r)}_0
+\widetilde{\resto} ^{(r)}_1, \\& \widetilde{\resto} ^{(r)}_0:= \sum _{\substack{|\mu |=|\nu | = r+1\\   |m|\le r} } a_{m\mu \nu
}^{(r)}
e^{\im m t}z^\mu \overline{z}^\nu  \\&  \widetilde{\resto} ^{(r)}_1:= \sum _{\substack{|\mu
|=|\nu |-1  = r\\   |m|\le r}  }
e^{\im m t} z^\mu \overline{z}^\nu \langle    \Phi
_{m\mu \nu }^{(r)}, f\rangle  + \sum _{\substack{|\nu
|=|\mu |-1  = r
\\   |m|\le r}    } e^{\im m t} z^\mu \overline{z}^\nu \langle \Psi _{m\mu \nu }^{(r)},
\overline{f}\rangle,
\end{aligned}
\end{equation}
i.e.  $\widetilde{\resto} ^{(r)}_0$ (resp. $
 \widetilde{\resto} ^{(r)}_1$)   defined as the sum of the terms in \eqref{eq:resto0} with $|\mu |=r$,  (resp. terms in \eqref{eq:resto1}
 with  $|\mu |+|\nu | =2r+1$).
Split $\widetilde{K}_{r+1}= K_{r+1}+ Z_{r+1}$,
   collecting inside $Z_{r+1}$ all the terms of $\widetilde{K}_{r+1}$
in the form either \eqref{eq:Z0}--\eqref{eq:Z01} or
\eqref{eq:Z1}--\eqref{eq:Z12}.  Apply Lemma
  \ref{sol.homo} with $\chi_{r+1}$ defined from $K_{r+1}$
  in the way \eqref{solhomo} is defined from \eqref{eq.stru.1}.
Then,

\begin{equation}
\label{HomEqmain} \left\{ H_F,\chi_{r+1}\right\} =-K_{r+1} .
\end{equation}
   $\chi_{r+1}$ is of form \eqref{chi.1} with $m_0=M_0=r$.
We can apply
  Lemmas \ref{lie_trans}--\ref{lem:lie_trans1} to $\chi_{r+1}$.
  Let $\mathfrak{F}_{r+1}$   be as $\mathfrak{F} $
  in \eqref{eq:Lie}.
For $ \Tr_{r+1}=\Tr_r\circ \mathfrak{F}_{r+1}$ set
\begin{equation}
\label{Hamr} H^{(r+1)}:=H^{(r)}\circ\mathfrak{F}_{r+1}=
H\circ(\Tr_r\circ\mathfrak{F}_{r+1}) = H\circ \Tr_{r+1}.
\end{equation}
Split
\begin{eqnarray}
\label{n1.1} H^{(r)}\circ\Tra_{r+1}&=& H_F+Z^{(r)} + Z_{r+1}
\\
\label{n1.3}
&+& (Z^{(r)}\circ\mathfrak{F}_{r+1}-Z^{(r)})\\
\label{n1.6} &+&
 {K}_{r+1}\circ\mathfrak{F}_{r+1}- {K}_{r+1}
\\
\label{n1.4} &+& H_F\circ\mathfrak{F}_{r+1}- \left(H_F+\left\{H_F,
\chi_{r+1} \right\}\right)
\\
\label{n1.5} &+& (\resto^{(r)}_{0 }- {\widetilde{\resto}^{(r)} _{0
}} +\resto^{(r)}_{1}- {\widetilde{\resto}^{(r)}_{1
}})\circ\mathfrak{F}_{r+1}
\\
\label{n1.7} &+&   (\resto
^{(r)}_2+...+\resto
^{(r)}_5)\circ\mathfrak{F}_{r+1}
\\ \label{n1.8}  &+&
  \resto ^{(r)}_6\circ\mathfrak{F}_{r+1}\\ \label{n1.9}  &+&
  \resto ^{(r)}_7\circ\mathfrak{F}_{r+1}
 \ .
\end{eqnarray}
Define $Z^{(r+1)}:=Z^{(r)}+Z_{r+1}$.  Then it is a degree $2r+2$
polynomial in $(z,\overline{z},f,\overline{f}) $ of the form
\eqref{e.12} satisfying \eqref{eq:Z0}--\eqref{eq:nf real}.  $ \resto ^{(r)}_7\circ\mathfrak{F}_{r+1}$ can be absorbed in $ \resto ^{(r+1)}_7 $.
For \eqref{n1.3}-- \eqref{n1.8} we apply Lemma
\ref{lem:lie_trans1} with $n_0+1= N+2$. The reminder terms corresponding to
$\mathbf{G}$ in \eqref{eq:LieTransf g} can be absorbed in $ \resto ^{(r+1)}_7 $.
Ti  all the terms corresponding to  $\mathcal{V}$
in \eqref{eq:LieTransf g} we can apply  (b) Lemma \ref{lem:lie_trans1}. In
particular, they are of the form \eqref{eq: LieTransf g1} with $L'=L+lr$
and with $|m'|\le L+ lr $. Since $L\ge 1$, this means they are of degree at least
$ 2r+ 4$ and they can be absorbed in $\resto ^{(r+1)}$.
\qed

\subsection{Application of Theorem \ref{theor:normal form}}
\label{subsec:appl nf}

We apply Theorem \ref{theor:normal form} for $r= N+1 $. Hence
we have $H^{(r)}:=H\circ \Tr_{r }=
H_F+Z_0^{(r)}+Z_1^{(r)}+\resto^{(r)}$. To simplify notation, but not
only for this reason, we will drop the super indexes but not in an
obvious way. The first step is obvious.  We set \begin{equation}
\label{eq:zr0}Z_0=Z^{(r)}_0.\end{equation} To proceed,  we go back
to the normal forms in (N1) and to conditions
\eqref{eq:Z1}--\eqref{eq:nf real}.
Set \begin{eqnarray}    &\mathbf{M}=
\{ (m, \mu
,\nu ) :  \lambda \cdot(\mu-\nu)-m  <-\underline{c} \, , \,  |\mu | =|\nu
|-1\, , \, |m |\le |\mu
 | \le N  \} \nonumber \\&  \label{eq:BfM}  \mathbf{M}'=
\{ (m',
\mu ',\nu ' ) :  (-m',
\nu ',\mu ' ) \in \mathbf{M}  \} .
\end{eqnarray}
Notice that the  two inequalities  $|\mu | =|\nu
|-1$ and   $ |m |\le |\mu
 | \le N $,  by  hypothesis (H7) imply  $ \lambda \cdot(\nu-\mu)+m  \neq \underline{c} $.

We have
\begin{equation}\label{eq:Z1r2}\begin{aligned}&  Z^{(r)}_1=
\sum _{  (m, \mu
,\nu)\in \mathbf{M}  } e^{\im   mt} z^\mu \overline{z} ^\nu\langle \Phi ^{(r)}_{m \mu
\nu}, f\rangle +\sum _{  (m', \mu '
,\nu ')\in \mathbf{M}'  } e^{\im   m't} z^{\mu'}\overline{z}^{\nu'}\langle \Psi
^{(r)} _{m ' \mu  ' \nu '}, \overline{f}\rangle .
\end{aligned}\nonumber
\end{equation}
By (iii.0)--(iii.1) Theorem \ref{theor:normal form}

\begin{equation}
\label{eq:real final} \overline{\Phi ^{(r)}} _{m\mu \nu }= \Psi ^{(r)}_{-m\nu \mu
} .
\end{equation}

\begin{definition}\label{def:M}
Denote by $M$ the subset of the indexes $(m,\mu , \nu )\in \mathbf{M}$ with the following
property: if $(m,\alpha , \beta )\in \mathbf{M}$   and if $\alpha _j\le \mu
_j$ and $\beta _j\le \nu _j$ for $j=0,...,n$, then $(\alpha , \beta
)= (\mu , \nu )$.
\end{definition}

\begin{definition}\label{def:Mprime}
Denote by $M'$ the subset of the indexes $(m',\mu ', \nu ')\in \mathbf{M}'$   with the
following property: if $(m',\alpha ', \beta ')\in \mathbf{M}'$   and if $\alpha
' _j\le \mu _j'$ and $\beta '_j\le \nu _j'$ for $j=0,...,n$, then
$(\alpha ', \beta ')= (\mu ', \nu ')$.
\end{definition}
The following is a straightforward consequence of  Theorem
\ref{theor:normal form}:
\begin{lemma}
\label{lemma:bjection} The following is a 1--1 and onto map $M\to
M'$:
\begin{equation}
\label{eq:isomorphism} (m, \mu , \nu )\to (m', \mu ' , \nu ') \text{
where } m'=-m, \, \mu '=\nu  , \, \nu '=\mu .
\end{equation}

\end{lemma}

We   set, and this is non obvious,
\begin{equation}
\label{eq:Z1fin}\begin{aligned}&  Z _1= \sum _{(m, \mu ,\nu )\in M }
e^{\im   mt} z^\mu \overline{z} ^\nu\langle \Phi ^{(r)}_{m \mu \nu},
f\rangle + \sum _{(m', \mu ',\nu ')\in M '}e^{\im   m't}
z^{\mu'}\overline{z}^{\nu'}\langle \Psi ^{(r)} _{m ' \mu  ' \nu '},
\overline{f}\rangle
\end{aligned}
\end{equation}
and
\begin{equation}
\label{eq:restofin}\begin{aligned}&  \resto = \resto ^{(r)} +
 Z^{(r)}_1-Z_1.
\end{aligned}
\end{equation}
Finally, to simplify notation, we set

\begin{equation}
\label{eq:coef fin} \Phi  _{m \mu \nu}=\Phi ^{(r)}_{m \mu \nu} \, ,
\quad \Psi  _{m \mu \nu}=\Psi ^{(r)}_{m \mu \nu}.
\end{equation}

 \section{Dispersion}
 \label{sec:dispersion}

We apply   Theorem \ref{theor:normal form} for $r= N+1 $ and we use
the notation in Subsection \ref{subsec:appl nf}.
 We
will show:
\begin{theorem}\label{proposition:mainbounds} There is a fixed
$C >0$ such that for $\varepsilon _0>0$ sufficiently small and for
$\epsilon \in (0, \varepsilon _0)$ we have
\begin{eqnarray}
&   \|  f \| _{L^r_t( [0,\infty ),W^{ 1 ,p}_x)}\le
  C \epsilon \text{ for all admissible pairs $(r,p)$}
  \label{Strichartzradiation}
\\& \| z ^{\mu +\nu} \| _{L^2_t([0,\infty ))}\le
  C \epsilon \text{ for all $(m, \mu , \nu )\in M$}   \label{L^2discrete}\\& \| z _j  \|
  _{W ^{1,\infty} _t  ([0,\infty )  )}\le
  C \epsilon \text{ for all   $j\in \{ 0, \dots , n\}$ }
  \label{L^inftydiscrete} .
\end{eqnarray}
\end{theorem}
By Theorem \ref{theor:normal form}, Theorem
\ref{proposition:mainbounds} implies Theorem \ref{theorem-1.1}.
Notice that \eqref{NLS1}  is time reversible, so in particular
\eqref{Strichartzradiation}--\eqref{L^inftydiscrete} are true over
the whole real line. The proof, though, exploits that $t\ge 0$,
specifically when for $\lambda \in \sigma _c(\mathcal{H})$  we
choose $R_{\mathcal{H}}^+(\lambda )=R_{\mathcal{H}} (\lambda +\im 0
)$ rather than $R_{\mathcal{H}}^-(\lambda )=R_{\mathcal{H}} (\lambda
-\im 0 )$ in formula \eqref{eq:g variable}. See the discussion on
p.18 \cite{SW3}.

The proof of
 Theorem
\ref{proposition:mainbounds}  is a  standard continuation
argument. We assume

\begin{eqnarray}
&   \|  f \| _{L^r_t([0,T],W^{ 1 ,p}_x)}\le
  C _1\epsilon \text{ for all admissible pairs $(r,p)$} \label{4.4a}
\\& \| z ^{\mu +\nu} \| _{L^2_t([0,T])}\le
 C_2 \epsilon  \text{ for all $(m, \mu , \nu )\in M$}  \label{4.4}\\& \| z _j  \|
  _{W ^{1,\infty} _t  ([0,T]  )}\le
  C_3 \epsilon \text{ for all   $j\in \{ 1, \dots , m\}$ } \label{4.4bis}
\end{eqnarray}
for fixed sufficiently large  constants $C_1$, $C_2$ and $C_3$. Then
  we prove that  for $\epsilon $ sufficiently small,
\eqref{4.4a}--\eqref{4.4bis} imply the same estimate but with $C_1$,
$C_2$, $C_3$ replaced by $C_1/2$, $C_2/2$, $C_3/2$. Then
\eqref{4.4a}--\eqref{4.4bis} hold with $[0,T]$ replaced by
$[0,\infty )$.

The proof consists in three main steps.
\begin{itemize}
\item[(i)] Estimate $f$ in terms of $z$.
\item[(ii)] Substitute the variable $f$  with a
new "smaller" variable $g$ and find smoothing estimates for $g$.
\item[(iii)] Reduce the system for $z$ to a closed system involving
only the $z$ variables, by insulating the  part of $f$  which
interacts with $z$, and by decoupling the rest (this reminder is
$g$). Then clarify the nonlinear Fermi golden rule.
\end{itemize}
These three steps are the same of the material in
\cite{bambusicuccagna} from Section 7 on,  and \cite{Cu1} from
Section 10 on.   We start by sketching  steps (i) and
(ii). Step (i) is encapsulated by the following proposition:

\begin{proposition}\label{Lemma:conditional4.2} Assume
\eqref{4.4a}--\eqref{4.4bis}. Then there exist constants $K_1$ and
$C=C(C_1,C_2,C_3)$,
  with $K_1$ independent of $C_1$, such that, if
  $C  \epsilon  $  is sufficiently small,   then we have
\begin{eqnarray}
&   \|  f \| _{L^r_t([0,T],W^{ 1 ,p}_x)}\le
  K_1  \epsilon \text{ for all admissible pairs $(r,p)$}\ .
  \label{4.5}
\end{eqnarray}
\end{proposition}
Consider $Z_1$ of the form \eqref{eq:Z1}--\eqref{eq:nf real}.   Then
we have (with finite sums)

\begin{equation}\label{eq:f variable} \begin{aligned}  &\im \dot f -
\mathcal{H}f   = \sum _{(m', \mu ', \nu ')\in M '} e^{\im tm'}
z^{\mu '}\overline{z}^{\nu '}
  \Psi _{m'\mu '\nu '} (x)   +
 \nabla _{\overline{f}} \resto .
\end{aligned}\end{equation}
The proof of Proposition \ref{Lemma:conditional4.2} is standard and
we skip it, see \cite{cuccagnamizumachi}. The dominating  term in the rhs of \eqref{eq:f variable}
is the   first line in the rhs. Notice also, that Theorem
\ref{proposition:mainbounds} implies   by standard arguments, see \cite{cuccagnamizumachi},

 \begin{equation}\label{scattering1}  \lim_{t\to +\infty}
\left \|   f (t) -
 e^{ \im t \Delta   \sigma_3}{f}_+  \right \|_{H^1}=0
 \end{equation}
for a  $  f_+\in H^1$ with $ \| {f}_+    \|_{H^1}\le C \epsilon$ and
for a real valued function $\theta \in C^1 (\mathbb{R},\mathbb{R})
.$

Step (ii) in the proof of Theorem \ref{proposition:mainbounds}
consists in introducing the variable

\begin{equation}
  \label{eq:g variable}
g=f+    \sum _{(m',\mu ' , \nu ') \in M'} e^{\im tm'} z^{\mu '}
\overline{z}^{\nu '}
  R ^{+}_{\mathcal{H}} (\lambda \cdot (\mu '-\nu ')-m')
  \Psi _{m'\mu ' \nu '} (x).
\end{equation}
Notice that by  Lemma 7 ch. XIII.8 \cite{reedsimon}, we have $g\in L^{2,-S}_x(\R ^3)$ for $S>1/2$.
Substituting the new variable  $g$ in \eqref{eq:f variable}, the
first line on the rhs of  \eqref{eq:f variable} cancels out. By an
easier version of Lemma 4.3 \cite{cuccagnamizumachi} we have:

\begin{lemma}\label{lemma:bound g} For $\epsilon$  sufficiently small and for and  $S>0$  sufficiently large, there exists
 $C_0=C_0(\mathcal{H},S)$  a fixed constant such that
\begin{equation} \label{bound:auxiliary}\| g
\| _{L^2_tL^{2,-S}_x}\le C_0 \epsilon + O(\epsilon
^2).\end{equation}
\end{lemma}
We skip the proof, which is standard, see \cite{cuccagnamizumachi}.
As in \cite{BP2,SW3} and subsequent literature, the part of $f$
which couples nontrivially with $z$ comes from  the polynomial in
$z$  in the summation in  rhs\eqref{eq:g variable}. In some sense,
$g$ and $z$
    decouple.

\subsection{The Fermi golden rule}
\label{subsec:FGR}

We proceed like \cite{bambusicuccagna,Cu1}, with in \eqref{eq:crunch}
  a different  Lyapunov   functional  than in  \cite{bambusicuccagna,Cu1}.
Set $R_{\mu \nu m}^\pm =R_{ \mathcal{H} }^\pm (\lambda\cdot (\mu
-\nu ) -m).$ System \eqref{eq:SystemH} with the new hamiltonian $
H^{(r)}$ yields (notice that  in the first line in \eqref{eq:FGR0}
summation is over $M$, while in the second line line in
\eqref{eq:FGR0}  is over $M'$)
\begin{equation}\label{eq:FGR0} \begin{aligned} &
\im \dot z _j-\lambda _j z _j
=
\partial _{\overline{z}_j}Z_0(z) +   \sum  _{(m, \mu , \nu )\in M}  \nu _je^{\im mt}\frac{z ^\mu
 \overline{ {z }}^ { {\nu} } }{\overline{z}_j} \langle f ,
  \Phi
_{m\mu \nu }\rangle
 \\ &  +  \sum  _{(m  , \mu  , \nu  )\in M'}  \nu _j e^{\im m t}
 \frac{z ^{\mu  }
 \overline{ {z }}^ { {\nu  } } }{\overline{z}_j }  \langle \overline{f} ,
  \Psi
_{m \mu   \nu  }\rangle      +
\partial _{  \overline{z} _j}  \resto   .
\end{aligned}  \end{equation}
We substitute \eqref{eq:g variable} in  \eqref{eq:FGR0}. Then we
rewrite \eqref{eq:FGR0}
\begin{eqnarray} \label{equation:FGR1}& \im \dot z _j-
\lambda _j z _j  =
\partial _{\overline{z}_j}Z_0(z) +   {\mathcal{E}}_j -
\\ &   \label{equation:FGR121}-\sum  _{\substack{
  (m, \mu , \nu )\in M\\ (m', \mu ', \nu ')\in M'}}  \nu _je^{\im (m+m')t}\frac{z  ^{\mu +\mu '}
 \overline{ {z }}^ { {\nu} +\nu ' } }{\overline{z}_j}
 \langle R^{+}_{\mu' \nu' m'} \Psi _{m'\mu' \nu'  } ,
  \Phi
_{m\mu \nu }\rangle \\ &  \label{equation:FGR13}- \sum  _{\substack{
  (m, \mu , \nu )\in M '\\ (m', \mu ', \nu ')\in M'}}  \nu _je^{\im (m-m')t}\frac{z  ^{\mu +\nu '}
 \overline{ {z }}^ { {\nu} +\mu ' } }{\overline{z}_j}
 \langle R^{-}_{\mu' \nu' m'} \overline{\Psi  _{m'\mu' \nu'  }} ,
  \Psi
_{ m\mu \nu } \rangle
\end{eqnarray}
with $  {\mathcal{E}}_j:= \text{rhs\eqref{eq:FGR0}}
-\text{\eqref{equation:FGR121}}- \text{\eqref{equation:FGR13}}
. $
 Most terms in \eqref{equation:FGR121} --\eqref{equation:FGR13}
    can
be eliminated through   new variables. Summing up only on the
subsets $ \mathcal{M}_1\subseteq M\times M'$ and $
\mathcal{M}_2\subseteq M'\times M'$, defined by the fact that   the denominators  in \eqref{equation:FGR2} are
non zero, we set
\begin{equation}\label{equation:FGR2}  \begin{aligned}   &
z _j =\zeta _j +\sum  _{\mathcal{M}_1} \frac{ \nu _je^{\im (m+m')t}z
^{\mu +\mu '}
 \overline{ {z }}^ { {\nu} +\nu ' } \langle R^{+}_{\mu' \nu' m'} \Psi _{m'\mu' \nu'  } ,
  \Phi
_{m\mu \nu }\rangle}{( m+m'-\lambda \cdot (\mu +\mu '-{\nu} -\nu
'))\overline{z}_j}
  \\&  +\sum  _{\mathcal{M}_2}  \frac{\nu _je^{\im (m-m')t}z  ^{\mu +\nu '}
 \overline{ {z }}^ { {\nu} +\mu ' }
 \langle R^{-}_{\mu' \nu' m'} \overline{\Psi  _{m'\mu' \nu'  }} ,
  \Psi
_{ m\mu \nu } \rangle  }{( m -m' -\lambda \cdot (\mu +\nu '-{\nu} -
\mu '))\overline{z}_j} .
  \end{aligned}
\end{equation}

\begin{lemma}
\label{lem:FGR2}
 In
\eqref{equation:FGR2} we have contributions from all terms in
\eqref{equation:FGR121} (resp. \eqref{equation:FGR13}) with $m\neq
-m' $  (resp. $m\neq m' $).  \end{lemma}
\proof  We consider only contributions from \eqref{equation:FGR121}. Contributions from \eqref{equation:FGR13} can be treated similarly.
It is enough to show that  for $m\neq
-m' $,    $(m,\mu , \nu )\in \mathbf{M}$ and  $(m ',\mu ', \nu ')\in \mathbf{M} '$, we have

\begin{equation}\label{eq:denomin} m+m'-\lambda \cdot (\mu +\mu '-{\nu} -\nu
')\neq 0  . \end{equation} By  the definition of  $\mathbf{M}$  and
$\mathbf{M}'$ in \eqref{eq:BfM}, we have $| \mu +\mu '-{\nu} -\nu '
|\le 4N+2$ and $|m+m'|\le 2N$. Then lhs\eqref{eq:denomin}$=0$  by
hypothesis (H8) would imply $m+m'=0$. This proves
\eqref{eq:denomin}. \qed

By \eqref{4.4}--\eqref{4.4bis} and by the fact that in \eqref{equation:FGR2}  we have $|
{\nu} |>1$, we obtain:

\begin{equation}  \label{equation:FGR3} \begin{aligned}   & \| \zeta  -
 z  \| _{L^2_t} \le \widetilde{C} C_3 \epsilon \sum _{ (m',\mu ', \nu ') \in M'
 } \| z ^{\mu ' +\nu ' }\| _{L^2_t}
\le CC_2C_3 \epsilon ^2,\\&  \| \zeta  -
 z \| _{L^\infty _t} \le C C_3^3\epsilon ^3 ,
\end{aligned}
\end{equation}
with $C$ a fixed constant. In the new variables, equation
\eqref{equation:FGR1}--\eqref{equation:FGR13}  is of the form ($ \mathcal{D} _j$ is discussed later, in the course of the proof of Lemma \ref{lemma:FGR1})

\begin{equation} \label{equation:FGR4} \begin{aligned} & \im \dot \zeta _j-
\lambda _j \zeta _j  -
\partial _{\overline{\zeta}_j}Z_0(\zeta ) -   \mathcal{D} _j =
\\ &     -\sum  _{\substack{ (m,\mu , \nu )\in M  \\
    (-m,\mu ', \nu ')\in M' \\  \lambda \cdot (\mu +\mu ')
    =
  \lambda \cdot (\nu +\nu ')}}  \nu _j \frac{\zeta  ^{\mu +\mu '}
 \overline{ {\zeta }}^ { {\nu} +\nu ' } }{\overline{\zeta}_j}
 \langle R^{+}_{\mu ' \nu ' (-m)} \overline{\Phi  _{  m
 \nu' \mu'  }} ,
  \Phi
_{m\mu \nu }\rangle  -\\& -\sum  _{\substack{ (m,\mu , \nu )\in M ' \\
    ( m,\mu ', \nu ')\in M' \\  \lambda \cdot (\mu -\nu  )
    =
  \lambda \cdot (\mu '- \nu ')}}  \nu _j \frac{\zeta  ^{\mu +\nu '}
 \overline{ {\zeta }}^ { {\nu} +\mu ' } }{\overline{\zeta}_j}
 \langle R^{-}_{\mu ' \nu ' m} \overline{\Psi  _{   m
 \mu' \nu'  }} ,
  \Psi
_{m\mu \nu }\rangle .\end{aligned}
\end{equation}
Set    $\alpha =\nu '$, $\beta =\mu '$  in the first term in
rhs\eqref{equation:FGR4}; replace in the second term in
rhs\eqref{equation:FGR4}, $\alpha =\nu '$, $\beta =\mu '$, $(\mu ,
\nu ) $ with $(\nu , \mu )$ and $m $ with $-m$. Then we get:

\begin{equation} \label{equation:FGR41} \begin{aligned} & \im \dot \zeta _j-
\lambda _j \zeta _j  -
\partial _{\overline{\zeta}_j}Z_0(\zeta ) -   \mathcal{D} _j =\\&
  -  \sum  _{\substack{ (m,\mu , \nu )\in M  \\
    ( m,\alpha , \beta )\in M  \\  \lambda \cdot (\mu - \nu )=
  \lambda \cdot ( \alpha - \beta )}}
  \nu _j \frac{\zeta  ^{\mu +\beta }
 \overline{ {\zeta }}^ { {\nu} +\alpha } }{\overline{\zeta}_j}
 \langle R^{+}_ \mathcal{H }(\lambda \cdot (\beta -\alpha) +m)
 \overline{\Phi  _{ m
 \alpha \beta  }} ,
  \Phi
_{m\mu \nu }\rangle -\\&  -
\sum  _{\substack{ (m,\mu , \nu )\in M  \\
    ( m,\alpha , \beta )\in M  \\  \lambda \cdot (\mu - \nu )=
  \lambda \cdot ( \alpha - \beta )}}
  \mu _j \frac{\zeta  ^{\nu +\alpha }
 \overline{ {\zeta }}^ { {\mu} +\beta } }{\overline{\zeta}_j}
 \langle R^{-}_ \mathcal{H }(\lambda \cdot (\beta -\alpha) +m)
  {\Phi  _{ m
 \alpha \beta  }} ,
  \overline{\Phi
_{m\mu \nu }}\rangle
       .
\end{aligned}
\end{equation}

Let now $X=\{ \lambda \cdot (\beta -\alpha) +m:  ( m,\alpha , \beta
)\in M  \} .$  $M$ is a finite set, so also $X$ is a finite set.
For each $w \in X$ let  $M_w= \{ ( m,\alpha , \beta
)\in M : \lambda \cdot (\beta -\alpha) +m=w \} .$

\begin{remark}
\label{rem:FGR4}  Notice that if $(m,\mu , \nu ) \in M _w$ and
$(m',\mu ', \nu ') \in M _w$, then $m=m'$ and $\lambda \cdot (\mu
-\nu )= \lambda \cdot (\mu '-\nu ')$  by (H8).
\end{remark}

 Recall Plemelji formula $\frac{1}{x\pm \im 0}
=P.V. \frac{1}{x}\mp \im \pi \delta (x)$. Then, proceeding as in \cite{SW3,bambusicuccagna} in this spot, we write

\begin{equation} \label{equation:FGR411} \begin{aligned} & \im \dot \zeta _j-
\lambda _j \zeta _j  -
\partial _{\overline{\zeta}_j}Z_0(\zeta ) -   \mathcal{D} _j =\\&
  -  \sum  _  {w\in X}   \sum  _{\substack{  (m,\mu , \nu )\in M _w \\
    ( m,\alpha , \beta )\in M  _w  }}
   \frac{\nu _j+\alpha _j }{\overline{\zeta}_j}
 \langle  P.V.\frac 1 {\mathcal{H}-w }
 \overline{\zeta  ^{\alpha}
 \overline{ {\zeta }}^ { \beta }\Phi  _{ m
 \alpha \beta  }} ,
  \zeta  ^{\mu}
 \overline{ {\zeta }}^ { \nu }\Phi
_{m\mu \nu }\rangle  \\& + \im \pi
 \sum  _ {w\in X}   \sum  _{\substack{ (m,\mu , \nu )\in M _w \\
    ( m,\alpha , \beta )\in M  _w }}
   \frac{\alpha _j-\nu _j }{\overline{\zeta}_j}
 \langle \delta (\mathcal{H}-w )
 \overline{\zeta  ^{\alpha}
 \overline{ {\zeta }}^ { \beta }\Phi  _{ m
 \alpha \beta  }} ,
  \zeta  ^{\mu}
 \overline{ {\zeta }}^ { \nu }\Phi
_{m\mu \nu }\rangle .
\end{aligned}
\end{equation}
Here we take a minor departure from \cite{bambusicuccagna}. We
prove:
\begin{lemma}
  \label{lem:Lyapunov}
Let \begin{equation} \label{equation:FGR43} \begin{aligned} &
\mathbf{\Phi}_{w}:= \sum  _{  (m,\mu , \nu )\in M _w   }
  \zeta  ^{\mu }
 \overline{ {\zeta }}^ { {\nu}   }  \Phi
_{m\mu \nu }
       .
\end{aligned}
\end{equation}
Then, we have  \begin{equation} \label{equation:FGR44}
\begin{aligned} & \sum _{j=0}^n\frac{1}{2} \partial _t |\zeta _j|^2
+\pi   \sum _{ w\in X }
    \langle \delta (\mathcal{H}-w)
  \overline{\mathbf{\Phi}_{w}}, \mathbf{\Phi}_{w}\rangle =
\sum _{j=0}^n
  \Im \left (  \mathcal{D} _j \overline{\zeta} _j\right )
.
\end{aligned}
\end{equation}
\end{lemma} \proof We multiply  \eqref{equation:FGR411}  by $ \overline{\zeta} _j$
and sum  on $j=0,...,n$ and take the imaginary part of the sum.   Then,
  \eqref{equation:FGR44}
is an immediate  consequence of  two cancelations,
\eqref{eq:cancel0} and \eqref{eq:cancel1} below,  and one identity,
\eqref{eq:cancel2} below. To finish the proof of Lemma \ref{lem:Lyapunov} we need to state and prove \eqref{eq:cancel0}, \eqref{eq:cancel1}   and
\eqref{eq:cancel2}.

 \noindent By
\eqref{eq:Z0}--\eqref{eq:Z01} and by $Z_0=\overline{Z_0}$  we get the  the first cancelation:

\begin{equation} \label{eq:cancel0}  \begin{aligned} & 2\im
\sum _{j=0}^n\Im (\overline{\zeta}_j\partial _{\overline{\zeta}_j}Z_0(\zeta )
) =\sum _{j=0}^n (\overline{\zeta}_j\partial
_{\overline{\zeta}_j}Z_0(\zeta )  - {\zeta}_j
 \partial _{ {\zeta}_j}Z_0(\zeta )   )\\& =\sum _{j=0}^n \sum _{|\mu |=| \nu |}
 (\mu _j-\nu _j) a_{\mu \nu}
 \zeta ^{\mu} \overline{\zeta} ^{\nu}  =\sum _{|\mu |=| \nu | }
 (|\mu |-|\nu |) a_{\mu \nu}
 \zeta ^{\mu} \overline{\zeta} ^{\nu} =0
       .
\end{aligned}
\end{equation}
We turn to the second cancelation. For $(m,\mu ,
\nu ) $ and $( m,\alpha , \beta )$ elements of $  M  _w\subset \mathbf{M}$,
by the definition of  $\mathbf{M}$ in \eqref{eq:BfM}
 we have
\begin{equation}\label{eq:cancel3} \sum _{j=0}^n (\mu _j -\nu _j)
=: |\mu |-|\nu|= -1  =|\alpha |-| \beta |=\sum _{j=0}^n (\alpha _j -\beta _j) .\end{equation} By \eqref{eq:cancel3} we get $\sum _{j=0}^n (\nu _j+\alpha _j ) = |\nu|+|\alpha |= |\beta |+|\mu | $. Hence

\begin{equation}\label{eq:cancel1}   \begin{aligned} &\Im
\sum  _{\substack{ w\in X\\ (m,\mu , \nu )\in M _w \\
    ( m,\alpha , \beta )\in M  _w  }}
   (|\nu|+|\alpha | )
 \langle  P.V.\frac 1 {\mathcal{H}-w }
 \overline{\zeta  ^{\alpha}
 \overline{ {\zeta }}^ { \beta }\Phi  _{ m
 \alpha \beta  }} ,
  \zeta  ^{\mu}
 \overline{ {\zeta }}^ { \nu }\Phi
_{m\mu \nu }\rangle =\\&  \Im
\sum  _{\substack{ w\in X\\ (m,\mu , \nu )\in M _w \\
    ( m,\alpha , \beta )\in M  _w  }}
    \frac{|\alpha |+|\beta  | +  |\mu |+|\nu | }{2}
 \langle  P.V.\frac 1 {\mathcal{H}-w }
 \overline{\zeta  ^{\alpha}
 \overline{ {\zeta }}^ { \beta }\Phi  _{ m
 \alpha \beta  }} ,
  \zeta  ^{\mu}
 \overline{ {\zeta }}^ { \nu }\Phi
_{m\mu \nu }\rangle   \\& =0.
\end{aligned}
\end{equation}
Having stated and  proved  the second cancelation \eqref{eq:cancel1}, we turn to the last identity, \eqref{eq:cancel3} below.
Rearranging, by  $ \Re \langle \delta (\mathcal{H}-w )
 \overline{v_1} ,
 v_2\rangle = \Re \langle \delta (\mathcal{H}-w )
 \overline{v_2} ,
 v_1\rangle$ and  by \eqref{eq:cancel3}, we get

\begin{equation} \label{eq:cancel2} \begin{aligned} &
\sum  _{\substack{ w\in X\\ (m,\mu , \nu )\in M _w \\
    ( m,\alpha , \beta )\in M  _w }}
   (|\alpha |-|\nu | )
\Re \langle \delta (\mathcal{H}-w )
 \overline{\zeta  ^{\alpha}
 \overline{ {\zeta }}^ { \beta }\Phi  _{ m
 \alpha \beta  }} ,
  \zeta  ^{\mu}
 \overline{ {\zeta }}^ { \nu }\Phi
_{m\mu \nu }\rangle =\\& \frac{1}{2}
\sum  _{\substack{ w\in X\\ (m,\mu , \nu )\in M _w \\
    ( m,\alpha , \beta )\in M  _w }}
   (|\alpha |-|\nu | +|\mu |-| \beta | )
\Re \langle \delta (\mathcal{H}-w )
 \overline{\zeta  ^{\alpha}
 \overline{ {\zeta }}^ { \beta }\Phi  _{ m
 \alpha \beta  }} ,
  \zeta  ^{\mu}
 \overline{ {\zeta }}^ { \nu }\Phi
_{m\mu \nu }\rangle \\& =-\sum  _{  w\in X } \langle \delta
(\mathcal{H}-w )\mathbf{\Phi}_{w}, \mathbf{\Phi}_{w} \rangle .
\end{aligned}
\end{equation}
 This yields Lemma \ref{lem:Lyapunov}. \qed
\begin{remark}
\label{rem:FGR} Formula \eqref{equation:FGR44} with Lemma \ref{lemma:FGR1} below
is the crucial structural result in the paper. This is a form of the so called {\it nonlinear Fermi golden rule}.
\end{remark}

\begin{lemma}
\label{lemma:FGR1} Assume  inequalities \eqref{4.4}--\eqref{4.4bis}.
Then for a fixed constant $c_0$ we have
\begin{eqnarray}\label{eq:FGR7}     \sum _{j=0}^n
\|\mathcal{D}_j \overline{\zeta} _j\|_{ L^1[0,T]}
 \le (1+C_2)c_0 \epsilon ^{2}
 . \end{eqnarray}
\end{lemma}
 We postpone the proof,  assume the conclusion and complete
 the proof of Theorem \ref{proposition:mainbounds}.
 We introduce now  the following key hypothesis.

\begin{itemize}
\item[(H9')]  We assume
that for some fixed constants for any vector $\zeta \in
\mathbb{C}^n$ we have:
\begin{equation}\label{eq:FGR}      \sum _{
w\in X }
    \langle \delta (\mathcal{H}-w)
  \overline{\mathbf{\Phi}_{w}}, \mathbf{\Phi}_{w}\rangle
 \approx \sum _{ (m,\mu , \nu )
   \in M}  | \zeta ^{\mu +\nu }  | ^2 .
\end{equation}
\end{itemize}
Now we complete the proof of Theorem \ref{proposition:mainbounds}.
Notice that  lhs\eqref{eq:FGR}$\ge 0.$ By
\eqref{equation:FGR44}--\eqref{eq:FGR}

\begin{equation} \label{eq:crunch}\begin{aligned}&  \sum _j   | z
_j(t)|^2 +\sum _{ (m,\mu , \nu )
   \in M}  \| z^{\mu +\nu } \| _{L^2(0,t)}^2\lesssim
\epsilon ^2+ C_2\epsilon ^2.\end{aligned}
\end{equation}
By \eqref{equation:FGR3} this implies $\| z ^{\mu +\nu } \|
_{L^2(0,t)}^2\lesssim \epsilon ^2+ C_2\epsilon ^2$ for all the above
multi indexes. So, from  $\| z ^{\mu +\nu } \| _{L^2(0,t)}^2\lesssim
  C_2^2\epsilon ^2$ we conclude $\| z ^{\mu +\nu } \|
_{L^2(0,t)}^2\lesssim  C_2\epsilon ^2$. This means that we can take
$C_2\approx 1$. This yields Theorem \ref{proposition:mainbounds}.

\bigskip

 {\bf Proof of Lemma \ref{lemma:FGR1}}. We have schematically
\begin{equation}
 \label{equation:FGR10}\begin{aligned} &
\mathcal{D}_j =\mathcal{E}_j+\overline{\partial} _jZ_0(t,z ) -\overline{\partial} _jZ_0(t,\zeta ) \\& +\sum _{k} \partial _{z_k}   \partial _{\overline{z}_j} \left ( \sum _{   \mathcal{M}_1}   z
^{\mu +\mu '}\overline{ {z }}^ { {\nu} +\nu ' } +  \sum  _{\mathcal{M}_2}   z  ^{\mu +\nu '}
 \overline{ {z }}^ { {\nu} +\mu ' } \right ) \text{
rhs\eqref{eq:FGR0}}_{k}
  \\& -\sum _{k} \partial _{\overline{z}_k}   \partial _{\overline{z}_j} \left ( \sum _{   \mathcal{M}_1}   z
^{\mu +\mu '}\overline{ {z }}^ { {\nu} +\nu ' } +  \sum  _{\mathcal{M}_2}   z  ^{\mu +\nu '}
 \overline{ {z }}^ { {\nu} +\mu ' } \right ) \overline{\text{
 rhs\eqref{eq:FGR0}}}
  _{k}
\end{aligned}
\end{equation}
where the exponents in the second line are the sums of the exponents in
\eqref{equation:FGR2}, where  $\text{
rhs\eqref{eq:FGR0}}_{k}$ is  just \eqref{eq:FGR0} when $j=k$
 and where \begin{equation}
\label{equation:FGR101}\begin{aligned}  &  \mathcal{E}_j=
\partial _{  \overline{z} _j}  \resto +\\& \sum  _{(m, \mu , \nu )\in M}  \nu _je^{\im mt}\frac{z ^\mu
 \overline{ {z }}^ { {\nu} } }{\overline{z}_j} \langle g ,
  \Phi
_{m\mu \nu }\rangle
  +  \sum  _{(m  , \mu  , \nu  )\in M'}  \nu _j e^{\im m t}
 \frac{z ^{\mu  }
 \overline{ {z }}^ { {\nu  } } }{\overline{z}_j }  \langle \overline{g} ,
  \Psi
_{m \mu   \nu  }\rangle   \\& -\sum  _{\substack{
  (m, \mu , \nu ) \not  \in   M \\ \lambda \cdot (\mu -\nu )-m<\underline{c}  \\
     |\mu |=  |\nu |-1\\ (m', \mu ', \nu ')\in M'}}  \nu _je^{\im (m+m')t}\frac{z  ^{\mu +\mu '}
 \overline{ {z }}^ { {\nu} +\nu ' } }{\overline{z}_j}
 \langle R^{+}_{\mu' \nu' m'} \Psi _{m'\mu' \nu'  } ,
  \Phi
_{m\mu \nu }\rangle \\ &   - \sum  _{\substack{
  (m, \mu , \nu )  \not  \in   M ' \\ \lambda \cdot (\mu -\nu )-m>\underline{c}  \\
     |\mu |=  |\nu |+1  \\ (m', \mu ', \nu ')\in M'}}  \nu _je^{\im (m-m')t}\frac{z  ^{\mu +\nu '}
 \overline{ {z }}^ { {\nu} +\mu ' } }{\overline{z}_j}
 \langle R^{-}_{\mu' \nu' m'} \overline{\Psi  _{m'\mu' \nu'  }} ,
  \Psi
_{ m\mu \nu } \rangle   .
\end{aligned}
\end{equation}

We now
estimate one by one the  terms in  \eqref{equation:FGR10}.
Lemma \ref{lemma:FGR1} is an immediate  consequence of Lemmas  \ref{lemma:E}-- \ref{lemma:third} below.
 \begin{lemma}\label{lemma:E} There are fixed $C_0$  and $ \epsilon _0>0$ such that
 for $\epsilon \in (0, \epsilon _0)$ we have
 \begin{equation}
  \begin{aligned}  & \| \mathcal{E}_j \overline{\zeta} _j\| _{L^1_t[0,T]}\le (1+C_2) C_0 \epsilon ^2, \\& \|   \mathcal{E}_j  \| _{L^2_t[0,T]}\le C_0 \epsilon  .
\end{aligned}
\end{equation}
\end{lemma}
\proof We have for  $C=C(C_1,C_2, C_3)$ \begin{equation}
 \label{equation:FGR102}\begin{aligned} &\| \overline{\zeta} _j\partial _{\bar
z _j}\resto    \| _{L^1_t}\le \| \partial _{   \overline{z}
_j}\resto \| _{L^1_t} \|   \overline{\zeta} _j\| _{L^\infty _t}\le
2C_3\epsilon\sum _{d=0}^{7} \| \partial _{   \overline{z} _j}\resto
_d \| _{L^1_t} \le C \epsilon ^3  \end{aligned}\end{equation} by
\eqref{4.4a}--\eqref{4.4bis}, \eqref{equation:FGR3}, by Theorem
\ref{theor:normal form} for $r=N+1$.

 We prove now that, for $\Phi \in \mathcal{S}(\R ^3, \C )$,
 we have
 \begin{equation}  \label{equation:FGR1031}\begin{aligned} &\sum  _{(m, \mu , \nu )\in M}  \nu _j \| e^{\im m t}
\frac{z ^{\mu}   z  ^{ {\nu} }}{\overline{z}_j}   \overline{\zeta} _j\langle \Phi , g
\rangle \| _{L^1_t}+\\&   \sum  _{(m  , \mu  , \nu  )\in M'}  \nu _j\| e^{\im m t}
\frac{z ^{\mu}   z  ^{ {\nu} }}{\overline{z}_j}   \overline{\zeta} _j\langle \Phi , \overline{g}
\rangle \| _{L^1_t} \le (1+C_2) C_0 \epsilon^2 \end{aligned}
 \end{equation}
 for a fixed $C_0$. These are the terms responsible for $C_2 c_0 \epsilon ^2$ in \eqref{eq:FGR7}. All the other terms can be incorporated in $  c_0 \epsilon ^2$.
 We focus on
\begin{equation}  \label{equation:FGR103} e^{\im m t}
\frac{z ^{\mu}   z  ^{ {\nu} }}{\overline{z}_j}   \overline{\zeta} _j\langle \Phi , g
\rangle =e^{\im m t}z ^{\mu}   z  ^{ {\nu} }  \langle \Phi , g \rangle +e^{\im m t}
\frac{z ^{\mu}   z  ^{ {\nu} }}{\overline{z}_j}  ( \overline{\zeta} _j - \overline{z} _j
)\langle \Phi , g \rangle . \end{equation}
We have  for   $ (m, \mu , \nu ) $ either in $ M$ or in $M'$ and by  Lemma \ref{lemma:bound g}

 \begin{equation}  \label{equation:FGR104} \| e^{\im m t}z ^{\mu}  \overline{z}  ^{ \nu }  \langle \Phi , g \rangle \|
_{L^1_t}\lesssim \| z ^{\mu}  \overline{z}  ^{ \nu }   \| _{L^2_t} \|
  g   \| _{L^2_tL_x^{2,-s}}\lesssim C_2 \epsilon ^2.\end{equation}
  Turning to the second term in rhs\eqref{equation:FGR103}
\begin{equation}  \label{equation:FGR104}\|  z^\mu \frac{\overline{z}  ^{ \nu }}{\overline{z}_j}  ( \overline{\zeta} _j -\overline{z}
_j )\langle \Phi , g \rangle \| _{L^1_t}\le  \|  z ^{\mu}  \frac{\overline{z}  ^{ \nu }}{\overline{z}_j}\| _{L^\infty _t}  \| \overline{\zeta} _j -\overline{z}
_j )  \| _{L^2_t}  \|
  g   \| _{L^2_tL_x^{2,-s}} \le   C (C_2 ,C_3)   \epsilon ^3 \nonumber \end{equation} by \eqref{4.4a}--\eqref{4.4bis},
by Lemma \ref{lemma:bound g}  and by
\eqref{equation:FGR3}.

We consider now contributions to $\mathcal{E}_j \overline{\zeta}_j $ coming from the third line of \eqref{equation:FGR101}. Terms  from the fourth line of \eqref{equation:FGR101}  can be treated similarly. We can write them

\begin{equation}  \label{equation:FGR105} \begin{aligned} & \sum  _{\substack{
  (m, \mu , \nu )  \not \in M \\ \lambda \cdot (\mu -\nu )  -m<c\\
     |\mu |=  |\nu | -1 \\ (m', \mu ', \nu ')\in M'}}  \nu _j ( \|   {z  ^{\mu +\mu '}
 \overline{ {z }}^ { {\nu} +\nu ' } }  \| _{L^1_t}  +
  \| \frac{z  ^{\mu +\mu '}
 \overline{ {z }}^ { {\nu} +\nu ' } }{\overline{z}_j}  (\overline{z}_j-\overline{\zeta} _j)\| _{L^1_t} ) \end{aligned}
\end{equation}
where we omitted the factors $e^{\im (m+m')t}$. It is easy to understand that the largest terms in \eqref{equation:FGR105} are the ones with $(m, \mu , \nu ) \in \mathbf{M} \backslash M $.  For each $(m, \mu , \nu ) \in \mathbf{M} \backslash M $  there exists a  $(m, \alpha , \beta )  \in M$ with $\alpha _j\le \mu
_j$ and $\beta _j\le \nu _j$ for all $j=0,...,n$ and with at least one of these not an equality. Hence

\begin{equation}  \label{equation:FGR106}\begin{aligned} &
\|   {z  ^{\mu +\mu '}
 \overline{ {z }}^ { {\nu} +\nu ' } }  \| _{L^1_t}\le \| z\| _{L^\infty _t}
 \| z^{\alpha +\beta}\| _{L^2_t} \| z^{\mu '+\nu '}\| _{L^2_t}
 \le C_3C_2^2 \epsilon ^3.\end{aligned}
\end{equation}

The second terms in \eqref{equation:FGR105}  can be bounded by
\begin{equation}  \label{equation:FGR107}
\|   \frac{z  ^{\mu +\mu '}
 \overline{ {z }}^ { {\nu} +\nu ' } }{\overline{z}_j}  (\overline{z}_j-\overline{\zeta} _j)  \| _{L^1_t}\le  \| \frac{z  ^{\mu  }
 \overline{ {z }}^ { {\nu}  } }{\overline{z}_j} \| _{L^\infty _t}
  \| z^{\alpha +\beta}\| _{L^2_t} \| z-\zeta \| _{L^2_t}
 \le C(C_2, C_3) \epsilon ^3.
\end{equation}
We have $\| \mathcal{E}_k  \|  _{L^2_t}\le C(C_2,C_3) \epsilon ^2 $,
as can be easily seen by \eqref{equation:FGR101} and the estimates used in
\eqref{equation:FGR102}.\qed

\begin{lemma}\label{lemma:chg variable} There is a fixed $\epsilon _0>0$ such that  assuming \eqref{4.4a}--\eqref{4.4bis} and for $\epsilon \in (0, \epsilon _0)$ we have

      \begin{equation}  \label{equation:FGR1071}\| (\overline{\partial} _jZ_0(t,z ) -\overline{\partial} _jZ_0(t,\zeta  ))  \zeta  _j\|
_{L^1_t[0,T]}\le  C(C_2, C_3) \epsilon ^3.\end{equation}
  \end{lemma}

  \proof
  We consider quantities  $(\frac{\zeta ^{\mu} \overline{\zeta}
^{\nu}}{\overline{\zeta}_j}-\frac{z ^{\mu} \overline{z}
^{\nu}}{\overline{z}_j}  )
 \overline{\zeta}  _j $ with $ (\mu , \nu )$ s.t. $|\mu |=| \nu |$
and $\lambda \cdot  (\mu - \nu ) =0, $ see \eqref{eq:Z01}.
  By
Taylor expansion these are
\begin{equation}\label{equation:FGR11} \sum \partial _k
(\frac{z ^{\mu} \overline{z}
^{\nu}}{\overline{z}_j}) (\zeta  _k-z _k)
\overline{\zeta} _j + \sum \overline{\partial} _k (\frac{z ^{\mu} \overline{z}
^{\nu}}{\overline{z}_j}) (\overline{\zeta} _k-\overline{z} _k)\overline{\zeta} _j +\overline{\zeta} _j  O(|z -\zeta |^2).
\end{equation}
It is straightforward that by \eqref{4.4bis} and \eqref{equation:FGR3} \begin{equation} \label{equation:FGR111}\|  \overline{\zeta} _j  O(|z -\zeta |^2) \| _{L^1_t} \le  C(C_2, C_3) \epsilon ^3.
\end{equation}
Turning to the first terms in \eqref{equation:FGR11} we split
\begin{equation}\label{equation:FGR112}\mu _k \frac{z ^{\mu} \overline{z}
^{\nu}}{z_k\overline{z}_j} \overline{\zeta} _j (\zeta  _k-z _k) =  \mu _k\frac {z ^{\mu} \overline{z}
^{\nu} }{z_k }  (\zeta  _k-z _k)  +  \mu _k \frac{z ^{\mu} \overline{z}
^{\nu}} {z_k\overline{z}_j} (\overline{\zeta} _j-z_j) (\zeta  _k-z _k) .
\end{equation}
The second term in rhs\eqref{equation:FGR112} can be treated like \eqref{equation:FGR111}. To bound the first  term in rhs\eqref{equation:FGR112}
we substitute \eqref{equation:FGR2}. We have
\begin{equation}
\label{equation:FGR113}\begin{aligned} & \mu _k\|
\frac{z ^{\mu} \overline{z}
^{\nu}}{z_k }   (\zeta  _k-z _k)\| _{L^1_t}
\le \sum  _{\mathcal{M}_1\cup \mathcal{M}_2 } \mu _k \beta  _k \|
\frac{z ^{\mu} \overline{z}
^{\nu}}{z_k }   e^{\im (m\pm m')t}\frac{z ^{\alpha + \alpha '} \overline{z}
^{\beta +\beta '}}{\overline{z}_k }\| _{L^1_t}.
\end{aligned}\nonumber
\end{equation}
We claim that  \eqref{eq:Z01}, $|\mu |+|\nu | \le 2 N+2 $, by (i)
 Theorem \ref{theor:normal form}, and
 hypothesis (H8)
   imply that there is at least one
index $\nu _\ell \neq 0$ such that $\lambda   _\ell =\lambda  _k$. Indeed, if
  $k \ge 1$,

\begin{equation}\label{equation:FGR1330}\lambda  _{k } \sum _{\ell : \lambda _\ell = \lambda  _{k }} (\mu  _{\ell}-\nu  _{\ell})=0
=\sum _{\ell : \lambda _\ell = \lambda  _{k }} (\mu  _{\ell}-\nu
_{\ell})  .
\end{equation}
This yields the existence of the $\nu   _\ell \neq 0$ for $k\neq 0$.
\eqref{equation:FGR1330} implies
\begin{equation}\label{equation:FGR13300}  \sum _{\ell \neq 0} (\mu _{\ell}-\nu _{\ell})=0 .
\end{equation}
Finally, \eqref{equation:FGR13300} and      $|\mu |=|\nu |$ imply also $ \mu _0 = \nu _0$, and our claim for $k=0$. Having established our claim, we can bound

\begin{equation}\label{equation:FGR12} \begin{aligned} &
\mu _k\beta  _k \|
\frac{z ^{\mu} \overline{z}
^{\nu}}{z_k  }    \frac{z ^{\alpha + \alpha '} \overline{z}
^{\beta +\beta '}}{\overline{z}_k }\| _{L^1_t} \le \mu _k\beta  _k \|
\frac{z ^{\mu} \overline{z}
^{\nu}}{z_k  \overline{z}_\ell}     \| _{L^\infty _t} \|
    \overline{z}_\ell \frac{z ^{\alpha  } \overline{z}
^{\beta }}{\overline{z}_k }\| _{L^2_t} \|
    {z ^{  \alpha '} \overline{z}
^{ \beta '}} \| _{L^2_t}
\\& \le  C(C_2, C_3)\epsilon ^4 \end{aligned}
\end{equation}
by the fact that the monomials in $Z_0$ have degree at least 4 (and so $ \|
\frac{z ^{\mu} \overline{z}
^{\nu}}{z_k  \overline{z}_\ell}     \| _{L^\infty _t}\le C (C_3) \epsilon ^2 $)
and that the monomial  $ \overline{z}_\ell \frac{z ^{\alpha  } \overline{z}
^{\beta }}{\overline{z}_k }$ belongs to the class of monomials with indexes in $M\cup M'$.
 Same
argument and bounds hold for   the second summation in
\eqref{equation:FGR11}. This yields \eqref{equation:FGR1071}. \qed

\begin{lemma}\label{lemma:third}There is a fixed $\epsilon _0>0$ such that  assuming \eqref{4.4a}--\eqref{4.4bis} and for $\epsilon \in (0, \epsilon _0)$ we have

      \begin{equation}   \|  \text{second}+\text{third line rhs\eqref{equation:FGR10}}\|
_{L^1_t[0,T]}\le  C(C_2, C_3) \epsilon ^3.\end{equation}

\end{lemma}
\proof
 We will only bound
\begin{equation}\label{equation:FGR130}\mu _k \nu _j\|
e^{\im (m+m')t}
\frac{z
^{\mu +\mu ' }  \overline{z} ^{  {\nu} +\nu '}}{z_k \overline{z}_j}\text{
rhs\eqref{eq:FGR0}}_{k} \overline{\zeta }_j\| _{L^1_t}
  \le  C(C_2, C_3)\epsilon ^3,
\end{equation}
 with  terms from the first line in rhs\eqref{equation:FGR2}.
 In particular
we assume $(m,\mu , \nu )\in M$ and   $(m'\mu ', \nu ')\in M'$. The other terms can be treated similarly.
To begin with, we will show

\begin{equation}\label{equation:FGR132} \mu _k \nu _j\| e^{\im (m+m')t}
\frac{z ^{\mu +\mu ' }  \overline{z} ^{  {\nu} +\nu '}}{z_k
\overline{z}_j} \overline{\zeta }_j {\partial} _{\overline{z}_k}
Z_0(t,z)  \| _{L^1_t}
  \le  C(C_2, C_3)\epsilon ^4.
\end{equation}
 It will be enough to consider

\begin{equation}\label{equation:FGR133}\mu _k \nu _j\nu ''_k\|
\frac{z
^{\mu +\mu ' }  \overline{z} ^{  {\nu} +\nu '}}{z_k \overline{z}_j}
  \frac{z
^{\mu '' }  \overline{z} ^{   \nu ''}}{  \overline{z}_k}  \overline{\zeta }_j\| _{L^1_t}
  \le  C(C_2, C_3)\epsilon ^4,
\end{equation}
with $(\mu '',\nu '')$ as in \eqref{eq:Z01}.  By the argument before
\eqref{equation:FGR13300} we can conclude that
\eqref{eq:Z01} and hypothesis (H8)  imply that there is at least one
index $\mu '' _\ell \neq 0$ such that $\lambda   _\ell =\lambda  _k$.
To prove  \eqref{equation:FGR133}
we substitute in the rhs  $\overline{\zeta }_j=\overline{z }_j+(\overline{\zeta }_j-\overline{z }_j).$ Then

\begin{equation}\label{equation:FGR134} \begin{aligned} &    \|
\frac{z
^{\mu +\mu ' }  \overline{z} ^{  {\nu} +\nu '}}{z_k \overline{z}_j}
  \frac{z
^{\mu '' }  \overline{z} ^{   \nu ''}}{  \overline{z}_k} (\overline{\zeta }_j-\overline{z }_j)\| _{L^1_t} \le \\&    \| \frac{z
^{\mu '' }  \overline{z} ^{   \nu ''}}{ z_\ell  \overline{z}_k} \| _{L^\infty _t}    \|
\frac{z
^{\mu +\mu ' }  \overline{z} ^{  {\nu} +\nu '} z_\ell }{z_k \overline{z}_j}
  \| _{L^2_t}   \| \zeta -z  \| _{L^2_t}
  \le  C(C_2, C_3)\epsilon ^5,\end{aligned}
\end{equation}
where we have used that  the first and the last factors
  in the second line of
\eqref{equation:FGR134} contribute at least  $\epsilon ^2$, and the middle one $\epsilon $. To complete \eqref{equation:FGR133} we prove
\begin{equation}\label{equation:FGR135}\mu _k \nu _j\nu ''_k\|
\frac{z
^{\mu +\mu ' }  \overline{z} ^{  {\nu} +\nu '}}{z_k  }
  \frac{z
^{\mu '' }  \overline{z} ^{   \nu ''}}{  \overline{z}_k}   \| _{L^1_t}
  \le  C(C_2, C_3)\epsilon ^4.
\end{equation}
We have by \eqref{4.4}--\eqref{4.4bis}
\begin{equation}\label{equation:FGR136} \begin{aligned} &    \|
\frac{z
^{\mu +\mu ' }  \overline{z} ^{  {\nu} +\nu '}}{z_k  }
  \frac{z
^{\mu '' }  \overline{z} ^{   \nu ''}}{  \overline{z}_k}\| _{L^1_t} \le
 \| \frac{z
^{\mu '' }  \overline{z} ^{   \nu ''}}{ z_\ell  \overline{z}_k} \| _{L^\infty _t}    \|
\frac{z
^{\mu +\mu ' }  \overline{z} ^{  {\nu} +\nu '} z_\ell }{z_k  }
  \| _{L^1_t}\le
\\&    C_3^2\epsilon ^2  \|
\frac{z
^{\mu   } z_\ell }{z_k  } \overline{z} ^{  {\nu}  }
  \| _{L^2_t}\|
 z
^{ \mu ' }  \overline{z} ^{  \nu '}
  \| _{L^2_t}
  \le   C_3^2 C_2^2   \epsilon ^4.\end{aligned}
\end{equation}
Hence we have proved \eqref{equation:FGR132}. Rewriting rhs$\eqref{eq:FGR0}=
\text{rhs}\eqref{equation:FGR1}+\eqref{equation:FGR121}+\eqref{equation:FGR13}
$, to complete the proof of \eqref{equation:FGR130} it is enough  to prove
\eqref{equation:FGR137}--\eqref{equation:FGR138} below.

We need  to show for $(m , \alpha , \beta  )$ either in $M$ or in $M'$,
  for  $(m' ,\alpha ', \beta ' ) $  in $M'$ and for $(\mu , \nu  ) $ sums of
  exponents in either of the two lines in \eqref{equation:FGR2},

\begin{equation}\label{equation:FGR137} \begin{aligned} & \mu _k \nu _j\beta _k\|
 \frac{z^\mu \overline{z}^\nu }{\overline{z}_j z_k}
  e^{\im (m  \pm m' )t}
  \frac{z  ^{\alpha  +\alpha '}\overline{ z}^ {\beta  +\beta ' }    }
  {\overline{z}_k  }  \overline{\zeta}_j  \|  _{L^1_t} \le  C(C_2, C_3)\epsilon ^4.
 \end{aligned}
\end{equation}
We also need to show:
\begin{equation}\label{equation:FGR138} \begin{aligned} &  \mu _k \nu _j\|
 \frac{z^\mu \overline{z}^\nu }{\overline{z}_j z_k} \mathcal{E}_k  \overline{\zeta}_j  \|  _{L^1_t} \le  C(C_2, C_3)\epsilon ^4.
 \end{aligned}
\end{equation}
Let us start with \eqref{equation:FGR137}. Substituting  $\overline{\zeta }_j=\overline{z }_j+(\overline{\zeta }_j-\overline{z }_j )$ and focusing for definiteness on  terms of the first line of \eqref{equation:FGR2},
we
reduce to

\begin{equation}\label{equation:FGR1371} \begin{aligned} & \mu _k \nu _j\beta _k\|
 \frac{z ^{\mu '+\mu ''} \overline{z}^{\nu '+\nu ''}  z  ^{\alpha  +\alpha '}\overline{ z}^ {\beta  +\beta ' }    }
  { z_k\overline{z}_k  }    \|  _{L^1_t} \le  C(C_2, C_3)\epsilon ^4 ,\\& \mu _k \nu _j\beta _k\|
 \frac{z ^{\mu '+\mu ''} \overline{z}^{\nu '+\nu ''}  z  ^{\alpha  +\alpha '}\overline{ z}^ {\beta  +\beta ' }    }
  {\overline{z}_j z_k\overline{z}_k  }  (z_j - \overline{\zeta}_j)  \|  _{L^1_t} \le  C(C_2, C_3)\epsilon ^4,
 \end{aligned}
\end{equation}
for some $(\widetilde{m}',\mu ' , \nu ')\in M$ and $(\widetilde{m}'',\mu '' , \nu '')\in M'$, with $ \mu _k$ (resp. $ \nu _k$) equal either to  $ \mu _k'$  or $ \mu _k''$ (resp. $ \nu _k'$  or $ \nu _k''$). Inequalities \eqref{equation:FGR1371} can be easily proved
using previous arguments. Finally let us turn now to  \eqref{equation:FGR138}.
Here too for definiteness we prove, for some $(\widetilde{m}',\mu ' , \nu ')\in M$ and $(\widetilde{m}'',\mu '' , \nu '')\in M'$,

\begin{equation}\label{equation:FGR1381} \begin{aligned} & \mu _k' \nu _j'\|
 \frac{z ^{\mu '+\mu ''} \overline{z}^{\nu '+\nu ''} }{\overline{z}_j z_k} \mathcal{E}_k  \overline{\zeta}_j  \|  _{L^1_t} \le  C(C_2, C_3)\epsilon ^4.
 \end{aligned}
\end{equation}
We have

\begin{equation}\label{equation:FGR1382} \begin{aligned} & \mu _k' \nu _j'\|
 \frac{z ^{\mu '+\mu ''} \overline{z}^{\nu '+\nu ''} }{\overline{z}_j z_k} \mathcal{E}_k  \overline{\zeta}_j  \|  _{L^1_t} \le \mu _k' \nu _j' \|
 \frac{z ^{\mu '+\mu ''} \overline{z}^{\nu '+\nu ''} }{\overline{z}_j z_k} \overline{\zeta}_j  \|  _{L^2_t}    \| \mathcal{E}_k  \|  _{L^2_t}  .
 \end{aligned}
\end{equation}
 We substitute $\overline{\zeta }_j=\overline{z }_j+(\overline{\zeta }_j-\overline{z }_j )$. Then we have

\begin{equation}\label{equation:FGR1382} \begin{aligned} & \mu _k'     \|
 \frac{z ^{\mu '+\mu ''} \overline{z}^{\nu '+\nu ''} }{  z_k} \|  _{L^2_t}
 \le   \mu _k'     \|
 \frac{z ^{\mu ' } \overline{z}^{\nu  '}  }{  z_k} \|  _{L^\infty_t}    \|
  z ^{ \mu ''} \overline{z}^{ \nu ''}  \|  _{L^2_t}
 \le C(C_2,C_3) \epsilon ^3,\\&  \mu _k' \nu _j '  \|
 \frac{z ^{\mu '+\mu ''} \overline{z}^{\nu '+\nu ''} }{\overline{z}_j z_k} (\overline{\zeta}_j -z_j) \|  _{L^2_t} \le  \mu _k' \nu _j ' \|
 \frac{z ^{\mu '+\mu ''} \overline{z}^{\nu '+\nu ''} }{\overline{z}_j z_k} \|  _{L^\infty _t} \| \overline{\zeta}  -z  \|  _{L^2_t}\\& \le C(C_2,C_3) \epsilon ^3.
 \end{aligned}
\end{equation}\qed

DISMI University of Modena and Reggio Emilia, Via Amendola 2, Pad.
Morselli, Reggio Emilia 42122, Italy.

{\it E-mail Address}: {\tt cuccagna.scipio@unimore.it}

\end{document}